\documentclass[12pt]{amsart}
\usepackage{latexsym, amsmath, amscd, amssymb, amsthm}
\usepackage{bm, mathrsfs, euscript}
\usepackage{psfrag, epsfig}
\newcommand{\scr}[1]{\bm{\EuScript{#1}}}
\newcommand{\ceil}[1]{\left\lceil{1}\right\rceil}
\DeclareMathOperator{\sat}{sat}
\DeclareMathOperator{\reach}{reach}
\DeclareMathOperator{\Ker}{Ker}
\DeclareMathOperator{\Aut}{Aut}
\DeclareMathOperator{\lcm}{lcm}
\DeclareMathOperator{\image}{Im}
\DeclareMathOperator{\reg}{reg}
\DeclareMathOperator{\pos}{pos}
\DeclareMathOperator{\supp}{supp}
\DeclareMathOperator{\negg}{neg}
\DeclareMathOperator{\initial}{in}
\DeclareMathOperator{\Spec}{Spec}
\DeclareMathOperator{\pdim}{pd}
\DeclareMathOperator{\interior}{int}
\DeclareMathOperator{\Pic}{Pic}
\newtheorem{lemma}{Lemma}[section] 
\newtheorem{theorem}[lemma]{Theorem}
\newtheorem{corollary}[lemma]{Corollary}
\newtheorem{proposition}[lemma]{Proposition}
\theoremstyle{definition}
\newtheorem{definition}[lemma]{Definition}
\newtheorem{example}[lemma]{Example}
\newtheorem{remark}[lemma]{Remark}

\newtheorem{problem}[lemma]{Problem}
\newtheorem*{acknowledgement}{Acknowledgements}
\numberwithin{equation}{section}
\renewcommand{\theequation}%
{\arabic{section}.\arabic{lemma}.\arabic{equation}}

\begin{document}

\title[Multigraded Regularity]{Multigraded Castelnuovo-Mumford
Regularity}

\author{Diane Maclagan} 
\address{Department of Mathematics \\ Stanford University \\ Stanford
\\ CA 94305 \\ USA}
\email{maclagan@math.stanford.edu}

\author{Gregory G. Smith}
\address{Department of Mathematics \\ Barnard College \\
  Columbia University \\ New York \\ \newline NY 10027 \\ USA}
\email{ggsmith@math.columbia.edu}


\begin{abstract}
We develop a multigraded variant of Castelnuovo-Mumford regularity.
Motivated by toric geometry, we work with modules over a polynomial
ring graded by a finitely generated abelian group. As in the standard
graded case, our definition of multigraded regularity involves the
vanishing of graded components of local cohomology.  We establish the
key properties of regularity: its connection with the minimal
generators of a module and its behavior in exact sequences.  For an
ideal sheaf on a simplicial toric variety $X$, we prove that its
multigraded regularity bounds the equations that cut out the
associated subvariety.  We also provide a criterion for testing if an
ample line bundle on $X$ gives a projectively normal embedding.
\end{abstract}

\maketitle

\section{Introduction}

Castelnuovo-Mumford regularity is a fundamental invariant in
commutative algebra and algebraic geometry.  Intuitively, it measures
the complexity of a module or sheaf.  The regularity of a module
approximates the largest degree of the minimal generators and the
regularity of a sheaf estimates the smallest twist for which the sheaf
is generated by its global sections.  Although the precise definition
may seem rather technical, its value in bounding the degree of
syzygies \cite{GLP} \cite{EL} and constructing Hilbert schemes
\cite{Gotzmann} \cite{iarrobinobook} has established that regularity
is an indispensable tool in both fields.

The goal of this paper is to develop a multigraded variant of
Castelnuovo-Mumford regularity.  We work with modules over a
polynomial ring graded by a finitely generated abelian group.
Imitating \cite{EG}, our definition of regularity involves the
vanishing of certain graded components of local cohomology.  We
establish the key properties of regularity: its connection with the
minimal generators of a module and its behavior in short exact
sequences.  As a consequence, we are able relate the regularity of a
module to chain complexes associated with the module.

Our primary motivation for studying regularity over multigraded
polynomial rings comes from toric geometry.  For a simplicial toric
variety $X$, the homogeneous coordinate ring, introduced in
\cite{Cox}, is a polynomial ring $S$ graded by the divisor class group
$G$ of $X$.  The dictionary linking the geometry of $X$ with the
theory of $G$-graded $S$\nobreakdash-modules leads to geometric
interpretations and applications for multigraded regularity.  We prove
that the regularity of an ideal sheaf bounds the multidegrees of the
equations that cut out the corresponding subvariety.  Multigraded
regularity also supplies a criterion for testing if an ample line
bundle on $X$ determines a projectively normal embedding.

Multigraded regularity consolidates a range of existing ideas.  In the
standard graded case, it reduces to Castelnuovo-Mumford regularity.
If $S$ has a nonstandard $\mathbb{Z}$\nobreakdash-grading, then our
definition is the version of regularity introduced in \cite{Benson}
for studying group cohomology.  When $S$ is the homogeneous coordinate
ring of a product of projective spaces, multigraded regularity is the
weak form of bigraded regularity defined in \cite{HoffmanWang}.  Our
description for the multigraded regularity of fat points
(Proposition~\ref{p:pointsReg}) is also connected with the results in
\cite{AdamTai}.  On the other hand, the versions of regularity
developed for Grassmannians in \cite{JayDeep} and abelian varieties in
\cite{PareschiPopa} do not fall under the umbrella of multigraded
regularity.

For ease of exposition we state our theorems in the case where $S$ is
the homogeneous coordinate ring of a smooth projective toric variety
$X$.  Let $B$ be the irrelevant ideal of $X$.  We denote by $\mathbb N
\scr{C}$ the semigroup generated by a finite subset $\scr{C} = \{
\bm{c}_{1}, \dotsc, \bm{c}_{\ell} \}$ of $G \cong \Pic(X)$.  In the
introduction, we restrict to the special case that $\scr{C}$ is the
minimal generating set of the semigroup of nef line bundles on $X$.
For example, if $X = \mathbb{P}^{d}$ then $G \cong \mathbb{Z}$, $S$
has the standard grading defined by $\deg(x_{i}) = \bm{1}$ for $1 \leq
i \leq n$, $B = \langle x_{1}, \dotsc, x_{n} \rangle$, and $\scr{C} =
\{ \bm{1} \}$.

The main point of this paper is to introduce the following definition,
which generalizes Castelnuovo-Mumford regularity.

\begin{definition} \label{i:def}
For $\bm{m} \in G$, we say that  a $G$-graded $S$-module $M$ is
$\bm{m}$-regular (with respect to $\scr{C}$) if the following
conditions are satisfied:
\begin{enumerate}
\item $H_{B}^{i}(M)_{\bm{p}} = 0$ for all $i \geq 1$ and all $\bm{p}
\in \bigcup ( \bm{m} - \lambda_{1} \bm{c}_{1} - \dotsb -
\lambda_{\ell} \bm{c}_{\ell} + \mathbb{N} \scr{C})$ where the union is
over all $\lambda_{1}, \dotsc, \lambda_{\ell} \in \mathbb{N}$ such
that $\lambda_{1} + \dotsb + \lambda_{\ell} = i-1$;
\item $H_{B}^{0}(M)_{\bm{p}} = 0$ for all $\bm{p} \in \bigcup_{1 \leq
j \leq \ell} ( \bm{m} + \bm{c}_{j} + \mathbb{N} \scr{C})$.
\end{enumerate}
We set $\reg(M) := \{ \bm{p} \in G : \text{$M$ is $\bm{p}$-regular}
\}$.
\end{definition}

In contrast with the usual notation, $\reg(M)$ is a set rather than a
single group element.  Traditionally, $G = \mathbb{Z}$ and the
regularity of $M$ refers to the smallest $\bm{m} \in G$ such that $M$
is $\bm{m}$-regular.  When $S$ has a multigrading, the group $G$ is
not equipped with a natural ordering so one cannot choose a smallest
degree $\bm{m}$.  More significantly, the set $\reg(M)$ may not even
be determined by a single element; see Example~\ref{ie:hir}.  From
this vantage point, bounding the regularity of a module $M$ is
equivalent to giving a subset of $\reg(M)$.

\begin{example} \label{ie:hir}
When $X$ is the Hirzebruch surface $\mathbb{F}_{2}$, the homogeneous
coordinate ring $S = \Bbbk[x_{1}, x_{2}, x_{3}, x_{4}]$ has the
$\mathbb{Z}^{2}$\nobreakdash-grading defined by $\deg(x_{1}) = \left[
\begin{smallmatrix} 1 \\ 0 \end{smallmatrix} 
\right]$, $\deg(x_{2}) = \left[
\begin{smallmatrix} - 2 \\ 0 \end{smallmatrix} 
\right]$, $\deg(x_{3}) = \left[
\begin{smallmatrix} 1 \\ 0 \end{smallmatrix} 
\right]$ and $\deg(x_{4}) = \left[ 
\begin{smallmatrix} 0 \\ 1 \end{smallmatrix} 
\right]$ and the irrelevant ideal $B$ is $\langle x_{1}x_{2},
x_{2}x_{3}, x_{3}x_{4}, x_{1}x_{4} \rangle$.  This grading identifies
 $\Pic(X)$ with $\mathbb{Z}^{2}$.  It follows that the semigroup of
nef line bundles is generated by the set $\scr{C} = \left\{ \left[
\begin{smallmatrix} 1 \\ 0 \end{smallmatrix} 
\right], \left[
\begin{smallmatrix} 0 \\ 1 \end{smallmatrix} 
\right] \right\}$.  A topological interpretation for $H_{B}^{i}(S)$
(see Section~3) shows that $\reg(S) = \bigl( \left[
\begin{smallmatrix} 1 \\ 0 \end{smallmatrix} 
\right] + \mathbb{N}^{2} \bigr) \cup \bigl( \left[
\begin{smallmatrix} 0 \\ 1 \end{smallmatrix} 
\right] + \mathbb{N}^{2} \bigr)$.  A graphical representation of
$\reg(S)$ appears in Figure~\ref{f:hirzreg}.  Observe that $\reg(S)$
is not determined by a single element of $G$ and $\left[
\begin{smallmatrix} 0 \\ 0 \end{smallmatrix} 
\right] \not\in \reg(S)$.
\end{example}

The following key result shows that the regularity of a module is a
cohomological approximation for the degrees of its minimal generators.

\begin{theorem} \label{i:gens}
Let $M$ be a finitely generated $G$-graded $S$-module.  If the module
$M$ is $\bm{m}$-regular, then the degrees of the minimal generators of
$M$ do not belong to the set $\reg(M) + \bigl( \bigcup_{1 \leq j \leq
\ell} (\bm{c}_{j} + \mathbb{N} \scr{C}) \bigr)$.
\end{theorem}

To emphasize the similarities between the multigraded form of
regularity and the original definition in \S14 of \cite{Mumford}, we
give the geometric translation of this result below.  For $\bm{p} \in
G$, $\mathscr{O}_{X}(\bm{p})$ is the associated line bundle on $X$.
Set $\mathscr{F}(\bm{p}) := \mathscr{F} \otimes
\mathscr{O}_{X}(\bm{p})$.

\begin{theorem} \label{i:sheaves}
Let $\mathscr{F}$ be a coherent $\mathscr{O}_{X}$-module.  If
$\mathscr{F}$ is $\bm{m}$-regular then for every $\bm{p} \in \bm{m} +
\mathbb{N} \scr{C}$ we have the following:
\begin{enumerate}
\item $\mathscr{F}$ is $\bm{p}$-regular;
\item the natural map $H^{0} \bigl( X, \mathscr{F}(\bm{p}) \bigr)
\otimes H^{0} \bigl( X, \mathscr{O}_{X}(\bm{q}) \bigr) \rightarrow
H^{0} \bigl( X, \mathscr{F}(\bm{p} + \bm{q}) \bigr)$ is surjective for
all $\bm{q} \in \mathbb{N} \scr{C}$;
\item $\mathscr{F}(\bm{p})$ is generated by its global sections.
\end{enumerate}
\end{theorem}

We highlight  two consequences of this theorem.  Firstly, if
$\mathscr{I}$ is an ideal sheaf on $X$ and $\bm{m} \in
\reg(\mathscr{I})$ then the subscheme of $X$ defined by $\mathscr{I}$
is cut out by equations of degree $\bm{m}$.  Secondly, if
$\mathscr{O}_{X}(\bm{p})$ is an ample line bundle and $\bm{p} \in
\reg(\mathscr{O}_{X})$ then the complete linear series associated to
$\mathscr{O}_{X}(\bm{p})$ gives a projectively normal embedding of
$X$.  In particular, if $\bm{0} \in \reg(\mathscr{O}_{X})$ then every
ample line bundle on $X$ gives a projectively normal embedding.

The next result illustrates a second key feature of regularity, namely
its behavior in exact sequences.  When $S$ has the standard grading,
the following are equivalent: the module $M$ is $\bm{m}$-regular, the
degrees of the $i$th syzgies are at most $\bm{m} + i$, and the
truncated module $M|_{\geq \bm{m}}$ has a linear resolution.  We
generalize these properties in the following way.

\begin{theorem} \label{i:complexes}
Let $d$ be the number of variable in polynomial ring $S$ minus the
rank of group $G$ and let $M$ be a finitely generated $G$-graded
$S$-module.
\begin{enumerate}
\item If $\;\; \dotsb \longrightarrow E_{3} \longrightarrow E_{2}
\longrightarrow E_{1} \longrightarrow E_{0} \xrightarrow{\;\;
\partial_{0} \;\;} M \longrightarrow 0$ is a chain complex of finitely
generated $G$-graded $S$-modules with $B$-torsion homology and
$\partial_{0}$ is surjective then
\[
\bigcup_{\phi \colon [d+1] \to [\ell]} \Bigl( \bigcap_{1 \leq i \leq
d+1} \bigl( - \bm{c}_{\phi(1)} - \dotsb - \bm{c}_{\phi(i)} +
\reg(E_{i}) \bigr) \Bigr) \subseteq \reg(M) \, .
\]
\item If $\bm{c} \in \reg(S) \cap \bigl( \bigcap_{1 \leq j \leq \ell}
( \bm{c}_{j} + \mathbb{N} \scr{C}) \bigr)$ and $\bm{m} \in \reg(M)$
then there exists a chain complex $\;\; \dotsb \longrightarrow E_{3}
\longrightarrow E_{2} \longrightarrow E_{1} \longrightarrow E_{0}
\longrightarrow M \longrightarrow 0$ with $B$-torsion homology and
$E_{i} = \bigoplus S(- \bm{m} - i \bm{c})$.
\end{enumerate}
\end{theorem}

If $S$ has the standard grading, then the inclusion in Part~1 is an
 equality when the $E_{i}$ form a minimal free resolution of $M$, and
 the chain complex in Part~2 is the minimal free resolution of
 $M|_{\geq \bm{m}}$.  Since $B$-torsion modules correspond to the zero
 sheaf on the corresponding toric variety, there is also a geometric
 version of this theorem involving regularity and resolutions of a
 sheaf.

The techniques used to establish Theorems~\ref{i:gens},
\ref{i:sheaves} and \ref{i:complexes} also apply to a larger class of
multigraded polynomial rings.  We develop these methods for pairs $(S,
B)$ where $S$ is a polynomial ring graded by a finitely generated
abelian group $G$, and  $B$ is a monomial ideal that encodes a certain
combinatorial structure of the grading.  This class of rings includes
the homogeneous coordinate rings of simplicial semi-projective toric
varieties.  We also establish these results for other choices of the
set $\scr{C}$.

This paper is organized as follows.  In the next section, we discuss
the basic definitions and examples of multigraded polynomial rings.
In Section~3, we establish some vanishing theorems for local
cohomology modules.  These results are based on a topological
description for the multigraded Hilbert series of $H_{B}^{i}(S)$.  The
definition of multigraded regularity is presented in Section~4.  We
also prove that in certain cases the definition of regularity is
equivalent to an apparently weaker vanishing condition.  Section~5
connects the multigraded regularity of a module with the degrees of
its minimal generators.  In Section~6, we reinterpret regularity of
$S$-module in terms of coherent sheaves on a simplicial toric variety
and study some geometric applications.  Finally, Section~7 examines
the relationship between chain complexes associated to a module or
sheaf and regularity.

\begin{acknowledgement}
We thank Dave Benson for conversations which expanded our notion of
regularity.  We are also grateful to Kristina Crona for conversations
about multigradings early in this project.  The computer software
package \texttt{Macaulay~2}~\cite{M2} was indispensable for computing
examples.  Both authors were partially supported by the Mathematical
Sciences Research Institute in Berkeley, CA.
\end{acknowledgement}

\section{Multigraded Polynomial Rings} \label{s:rings}

In this section, we develop the foundations of multigraded polynomial
rings.  Let $\Bbbk$ be a field and let $G$ be a finitely generated
abelian group.  Throughout this paper, we work with a pair $(S, B)$
where $S := \Bbbk[x_{1}, \dotsc, x_{n}]$ is a $G$-graded polynomial
ring and $B$ is a monomial ideal in $S$.  For a positive integer $m$,
we write $[m]$ for the set $\{ 1, \dotsc, m \}$.  The convex cone
generated by (or positive hull of) the vectors $\{\bm{v}_{1}, \dotsc,
\bm{v}_{m}\}$ is the set $\pos( \bm{v}_{1}, \dotsc, \bm{v}_{m} ) := \{
\lambda_{1} \bm{v}_{1} + \dotsb + \lambda_{m} \bm{v}_{m} :
\lambda_{1}, \dotsc, \lambda_{m} \in \mathbb{R}_{\geq 0} \}$.

A $G$-grading of the polynomial ring $S = \Bbbk[x_{1}, \dotsc, x_{n}]$
corresponds to a semigroup homomorphism $\mathbb{N}^{n}
\longrightarrow G$.  We say a monomial $\bm{x}^{\bm{u}}$ has degree
$\bm{p}$ if $\bm{u} \longmapsto \bm{p} \in G$.  This map induces a
decomposition $S = \bigoplus_{\bm{p} \in G} S_{\bm{p}}$ satisfying
$S_{\bm{p}} \cdot S_{\bm{q}} \subseteq S_{\bm{p} + \bm{q}}$ where
$S_{\bm{p}}$ is the $\Bbbk$-span of all $\bm{x}^{\bm{u}}$ of degree
$\bm{p}$.  The $G$-grading is determined by the set $\scr{A} := \{
\bm{a}_{1}, \dotsc, \bm{a}_{n} \}$ where $\bm{a}_{i} := \deg(x_{i})
\in G$ for $1 \leq i \leq n$.  We write $\mathbb{N} \scr{A}$ for the
subsemigroup of $G$ generated by $\scr{A}$.  Let $r$ be the rank of
$G$ and identify $\mathbb{R} \otimes_{\mathbb{Z}} G$ with
$\mathbb{R}^{r}$.  If $\overline{\bm{a}}_{i}$ denotes the image of
$\bm{a}_{i}$ in $\mathbb{R} \otimes_{\mathbb{Z}} G = \mathbb{R}^{r}$,
then the set $\overline{\scr{A}} := \{ \overline{\bm{a}}_{1}, \dotsc,
\overline{\bm{a}}_{n} \}$ is an integral vector configuration in
$\mathbb{R}^{r}$.  The monomial ideal $B$ corresponds to a chamber
(maximal cell) $\Gamma \subset \mathbb{R}^{r}$ in the chamber complex
of the vector configuration $\overline{\scr{A}}$.  The chamber complex
is the coarsest fan with support $\pos(\overline{\scr{A}})$ that
refines all the triangulations of $\overline{\scr{A}}$; see
\cite{BGS}.  We encode the choice of a chamber $\Gamma$ in the
monomial ideal 
\[
B := \bigl\langle \textstyle\prod\nolimits_{i \in \sigma} x_{i} :
\text{$\sigma \subseteq [n]$ with $\Gamma \subseteq \pos(
\overline{\bm{a}}_i: i \in \sigma )$} \bigr\rangle \, .
\]

Alternatively, the monomial ideal $B$ can be described by a regular
triangulation (Definition~5.3 in \cite{Ziegler}) of the dual vector
configuration.  If $d := n - r$, then the dual vector configuration is
a set of vectors $\overline{\scr{B}} := \{ \overline{\bm{b}}_{1},
\dotsc, \overline{\bm{b}}_{n} \}$ in $\mathbb{R}^{d}$ such that
\begin{equation} \label{dualconfiguration}
0 \longrightarrow \mathbb{R}^{d} \xrightarrow{\;\; [
\overline{\bm{b}}_{1} \; \dotsb \; \overline{\bm{b}}_{n}
]^{\textsf{T}} \;\;} \mathbb{R}^{n} \xrightarrow{\;\; [
\overline{\bm{a}}_{1} \; \dotsb \; \overline{\bm{a}}_{n} ] \;\;}
\mathbb{R}^{r} \longrightarrow 0
\end{equation}
is a short exact sequence.  The set $\overline{\scr{B}}$ is uniquely
determined up to a linear change of coordinates on $\mathbb{R}^{d}$;
see \S6.4 in \cite{Ziegler}.  We identify a triangulation of
$\overline{\scr{B}}$ with a simplicial complex $\Delta$.  Gale
duality, specifically Theorem~3.1 in \cite{BGS}, implies that the
chamber $\Gamma \subset \mathbb{R}^{r}$ corresponds to a regular
triangulation $\Delta$ of $\overline{\scr{B}}$ .  For $\sigma
\subseteq [n]$, let $\widehat{\sigma}$ denote the complement of
$\sigma$ in $[n]$. From this standpoint, we have $B = \bigl\langle
\prod\nolimits_{i \in \widehat{\sigma}} x_{i} : \text{$\sigma \in
\Delta$} \bigr\rangle$.

The simplicial complex $\Delta$ (or equivalently the chamber $\Gamma$)
also gives rise to two important subsemigroups of $\mathbb{N}
\scr{A}$.  The first subsemigroup $\scr{K}$ is the intersection
$\bigcap_{\sigma \in \Delta} \mathbb{N} \scr{A}_{\widehat{\sigma}}$
where $\scr{A}_{\widehat{\sigma}} := \{ \bm{a}_{i} : i \in
\widehat{\sigma} \}$.  Let $\mathbb{Z} \scr{A}$ be the subgroup of $G$
generated by $\scr{A}$.  The second subsemigroup $\scr{K}^{\sat}$ is
the saturation (or normalization) of $\scr{K}$ in $\mathbb{Z}
\scr{A}$.  In other words, $\scr{K}^{\sat}$ is the set of all $\bm{p}
\in \mathbb{N} \scr{A}$ such that the image $\overline{\bm{p}} \in
\mathbb{R}^{r}$ lies in the chamber $\Gamma$.  The interior of
$\scr{K}^{\sat}$, denoted $\interior \scr{K}^{\sat}$, consists of all
$\bm{p} \in \scr{K}^{\sat}$ such that $\overline{\bm{p}}$ lies in the
interior of $\Gamma$.

\begin{example} \label{e:standard}
Let $G = \mathbb{Z}$ and assume that $\bm{a}_{i} > 0$ for all $i$.
Since $G$ is torsion-free, we may identify $\bm{a}_{i}$ with
$\overline{\bm{a}}_{i}$.  The polynomial ring $S = \Bbbk[x_{1},
\dotsc, x_{n}]$ has the $\mathbb{Z}$-grading induced by $\deg(x_{i}) =
\bm{a}_{i}$ for $1 \leq i \leq n$.  The chamber complex of $\scr{A}$
has a unique maximal cell $\Gamma = \mathbb{R}_{\geq 0}$ and the
corresponding ideal is $B = \langle x_{1}, \dotsc, x_{n} \rangle$.
The dual vector configuration $\overline{\scr{B}}$ is given by the
columns of the matrix
\[
\begin{bmatrix}
\bm{a}_{n} & 0 & \dotsb & 0 & - \bm{a}_{1} \\
0 & \bm{a}_{n} & \dotsb & 0 & - \bm{a}_{2} \\
\vdots & \vdots & \ddots & \vdots & \vdots \\
0 & 0 & \dotsb & \bm{a}_{n} & - \bm{a}_{n-1} \\
\end{bmatrix} \, .
\]
The semigroup $\scr{K} \subset \mathbb{Z}$ consists of all nonnegative
multiples of $\lcm(\scr{A})$ and $\scr{K}^{\sat} = \mathbb{N}$. It
follows that $\scr{K}^{\sat} = \scr{K}$ if and only if $\bm{a}_{i} =
1$ for all $1 \leq i \leq n$.  When $\bm{a}_{i} = 1$ for all $i$, $S$
is standard graded polynomial ring and $B$ is the unique graded
maximal ideal. 
\end{example}

\begin{example} \label{e:hirzebruch}
Fix $t \in \mathbb{N}$.  Suppose that $G = \mathbb{Z}^{2}$ and let
$\scr{A}$ correspond to the columns of the matrix
\[
\begin{bmatrix} 
1 & - t & 1 & 0 \\ 0 & 1 & 0 &1
\end{bmatrix}  \, .
\] 
Again, we identify $\scr{A}$ with $\overline{\scr{A}}$.  The
polynomial ring $S = \Bbbk[x_{1}, x_{2}, x_{3}, x_{4}]$ has the
$\mathbb{Z}^{2}$\nobreakdash-grading defined by $\deg(x_{1}) = \left[
\begin{smallmatrix} 1 \\ 0 \end{smallmatrix} 
\right]$, $\deg(x_{2}) = \left[
\begin{smallmatrix} - t \\ 0 \end{smallmatrix} 
\right]$, $\deg(x_{3}) = \left[
\begin{smallmatrix} 1 \\ 0 \end{smallmatrix} 
\right]$ and $\deg(x_{4}) = \left[ 
\begin{smallmatrix} 0 \\ 1 \end{smallmatrix} 
\right]$.  There are two possible choices for the chamber $\Gamma$:
$\mathbb{R}_{\geq 0}^{2} = \pos \left( \left[ 
\begin{smallmatrix} 1 \\ 0 \end{smallmatrix} \right], \left[ 
\begin{smallmatrix} 0 \\ 1 \end{smallmatrix} \right] \right)$ 
or $\pos \left( \left[ 
\begin{smallmatrix} -t \\ 1 \end{smallmatrix} \right], \left[ 
\begin{smallmatrix} 0 \\ 1 \end{smallmatrix} \right] \right)$ 
and the ideal $B$ equals $\langle x_{1}x_{2}, x_{2}x_{3}, x_{3}x_{4},
x_{1}x_{4} \rangle$ or $\langle x_{1} x_{2}, x_{2}x_{3}, x_{2}x_{4}
\rangle$ respectively.  The dual vector configuration is given by the
columns of the matrix
\[
\begin{bmatrix} 
1 & 0 & -1 & 0 \\ 0 & 1 & t & -1
\end{bmatrix} \, .
\]  
Figure~\ref{f:Hirzfig} illustrates the associated vector
configurations when $t = 2$.  Regardless of the choice of $\Gamma$,
both $\scr{K}^{\sat}$ and $\scr{K}$ equal $\Gamma \cap
\mathbb{Z}^{2}$.
\begin{figure}[ht]
\psfrag{1}{\scriptsize $1$}
\psfrag{2}{\scriptsize $2$}
\psfrag{3}{\scriptsize $3$}
\psfrag{4}{\scriptsize $4$}
\psfrag{A}{$\overline{\scr{A}}$}
\psfrag{B}{$\overline{\scr{B}}$}
\epsfig{file=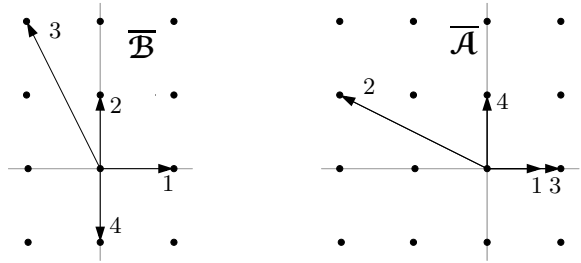, width=3in}
\caption{The vector configurations for Example~\ref{e:hirzebruch}
\label{f:Hirzfig}}
\end{figure}
\end{example}

\begin{example} \label{e:withtorsion}
If $G = \mathbb{Z}^{2} \oplus (\mathbb{Z}/2 \mathbb{Z})^{2}$ and
$\scr{A}$ corresponds to the columns of the matrix 
\[
\begin{bmatrix} 
1 & 0 & 1 & 0 & 1 \\ 
0 & 1 & 0 & 1 & 1 \\ 
1 & 0 & 1 & 0 & 0 \\ 
1 & 1 & 0 & 0 & 0
\end{bmatrix}
\] 
where the entries in bottom two rows are elements of $\mathbb{Z}/2
\mathbb{Z}$, then $S = \Bbbk[x_{1}, x_{2}, x_{3}, x_{4}, x_{5}]$ is
graded by a group with torsion.  The vector configuration
$\overline{\scr{A}}$ is given by the columns of the matrix 
$\left[
\begin{smallmatrix} 
1 & 0 & 1 & 0 & 1 \\ 0 & 1 & 0 & 1 & 1 
\end{smallmatrix} \right]$. 
The chamber $\Gamma$ is either $\Gamma_{1} := \pos \left( \left[
\begin{smallmatrix} 1 \\ 0 \end{smallmatrix} 
\right], \left[ 
\begin{smallmatrix} 1 \\ 1 \end{smallmatrix} 
\right] \right)$ or $\Gamma_{2} := \pos \left( \left[
\begin{smallmatrix} 0 \\ 1 \end{smallmatrix} 
\right], \left[ 
\begin{smallmatrix} 1 \\ 1 \end{smallmatrix} 
\right] \right)$ and the ideal $B$ equals either $\langle x_{1}x_{2},
x_{1}x_{4}, x_{1}x_{5}, x_{2}x_{3}, x_{3}x_{4}, x_{3}x_{5} \rangle$ in
the first case or $\langle x_{1}x_{2}, x_{2}x_{3}, x_{2}x_{5},
x_{1}x_{4}, x_{3}x_{4}, x_{4}x_{5} \rangle$ in the second.  The dual
vector configuration $\overline{\scr{B}}$ is given by the columns of
the matrix
\[
\begin{bmatrix}
1 & 1 & -1 & -1 & 0 \\
1 & -1 & -1 & 1 & 0 \\
1 & 1 & 1 & 1 & -2 
\end{bmatrix} \, .
\]
The two triangulations are illustrated in Figure~\ref{f:torsioneg}.
\begin{figure}[ht]
\psfrag{1}{\scriptsize $1$}
\psfrag{2}{\scriptsize $2$}
\psfrag{3}{\scriptsize $3$}
\psfrag{4}{\scriptsize $4$}
\psfrag{5}{\scriptsize $5$}
\psfrag{=}{{\scriptsize $5 = \infty$}}
\psfrag{G1}{$\Gamma_{1}$}
\psfrag{G2}{$\Gamma_{2}$}
\psfrag{A}{$\overline{\scr{A}}$}
\epsfig{file=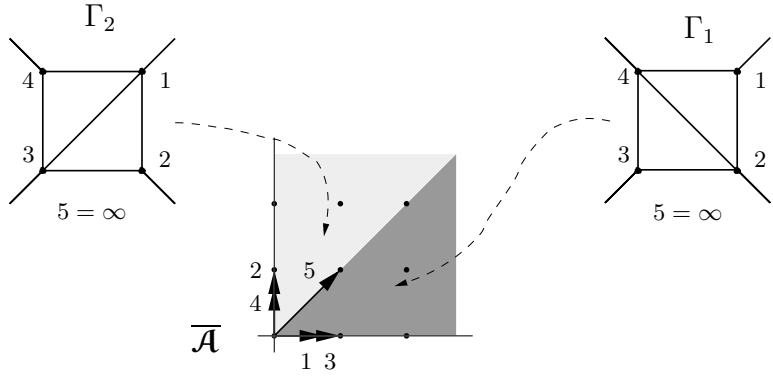, width=4in}
\caption{The two triangulations and chamber complex for
Example~\ref{e:withtorsion} \label{f:torsioneg}}
\end{figure}
If $\Gamma_{1}$ is the chosen chamber, then
\[
\scr{K} = \mathbb{N} \left\{ \left[
\begin{smallmatrix} 2 \\ 0 \\ 0 \\ 0 \end{smallmatrix} 
\right], \left[
\begin{smallmatrix} 2 \\ 2 \\ 0 \\ 0 \end{smallmatrix} 
\right] \right\} \qquad \text{and} \qquad 
\scr{K}^{\sat} = \mathbb{N} \left\{ \left[
\begin{smallmatrix} 1 \\ 0 \\ 0 \\ 0 \end{smallmatrix} 
\right], \left[
\begin{smallmatrix} 1 \\ 1 \\ 0 \\ 0 \end{smallmatrix} 
\right], \left[
\begin{smallmatrix} 0 \\ 0 \\ 1 \\ 0 \end{smallmatrix} 
\right], \left[
\begin{smallmatrix} 0 \\ 0 \\ 0 \\ 1 \end{smallmatrix} 
\right] \right\} \, .
\]
On the other hand, if $\Gamma_{2}$ is the chosen chamber, then
\[
\scr{K} = \mathbb{N} \left\{ \left[
\begin{smallmatrix} 0 \\ 2 \\ 0 \\ 0 \end{smallmatrix} 
\right], \left[
\begin{smallmatrix} 2 \\ 2 \\ 0 \\ 0 \end{smallmatrix} 
\right] \right\} \qquad \text{and} \qquad  
\scr{K}^{\sat} = \mathbb{N} \left\{ \left[
\begin{smallmatrix} 0 \\ 1 \\ 0 \\ 0 \end{smallmatrix} 
\right], \left[
\begin{smallmatrix} 1 \\ 1 \\ 0 \\ 0 \end{smallmatrix} 
\right], \left[
\begin{smallmatrix} 0 \\ 0 \\ 1 \\ 0 \end{smallmatrix} 
\right], \left[
\begin{smallmatrix} 0 \\ 0 \\ 0 \\ 1 \end{smallmatrix} 
\right] \right\} \, .
\]
\end{example}

The following lemma further illustrates the connections between the
ideal $B$ and the semigroups $\scr{K}$ and $\scr{K}^{\sat}$.  For
$\bm{p} \in G$, let $\langle S_{\bm{p}} \rangle$ denote the ideal
generated by all the monomials of degree $\bm{p}$ in $S$.  For $\bm{u}
\in \mathbb{N}^{n}$, let $\supp(\bm{u}) := \{ i : u_{i} \neq 0 \}
\subseteq [n]$.

\begin{lemma} \label{l:BandC}
If $\bm{p}$ belongs to the interior of $\scr{K}^{\sat}$ and
$\deg(\bm{x}^{\bm{u}}) = \bm{p}$ then $\bm{x}^{\bm{u}}$ belongs to the
ideal $B$.  Moreover, if $\bm{p} \in \scr{K}$, then $B \subseteq
\sqrt{\langle S_{\bm{p}} \rangle}$.
\end{lemma}

\begin{proof}
If $\bm{p} \in \interior \scr{K}^{\sat}$ and $\bm{x}^{\bm{u}} \in
S_{\bm{p}}$ then $\pos \bigl( \overline{\bm{a}}_{i} : i \in
\supp(\bm{u}) \bigr) \cap \interior \Gamma \neq \emptyset$.  Since
$\Gamma$ is a maximal cell in the chamber complex of
$\overline{\scr{A}}$, we have $\Gamma \subseteq \pos \bigl(
\overline{\bm{a}}_{i} : i \in \supp(\bm{u}) \bigr)$ and
$\dim_{\mathbb{R}} \pos \bigl( \overline{\bm{a}}_{i} : i \in
\supp(\bm{u}) \bigr) = r$.  Caratheodory's Theorem (Proposition~1.15
in \cite{Ziegler}) implies that $\pos\bigl( \overline{\bm{a}}_{i} : i
\in \supp(\bm{u}) \bigr)$ is the union of the $\pos(
\overline{\bm{a}}_{i} : i \in \sigma )$ where $|\sigma| = r$ and
$\sigma \subseteq \supp(\bm{u})$.  Hence, there exists $\sigma
\subseteq \supp(\bm{u})$ such that $|\sigma| = r$ and $\pos(
\overline{\bm{a}}_{i} : i \in \sigma) \cap \interior \Gamma \neq
\emptyset$.  Again, because $\Gamma$ is a chamber, we have $\Gamma
\subseteq \pos( \overline{\bm{a}}_{i} : i \in \sigma)$.  It follows
that the monomial $\prod_{i \in \sigma} x_{i}$ belongs to $B$ and
divides $\bm{x}^{\bm{u}}$.  This establishes the first assertion.

If $\bm{p} \in \scr{K}$, then for every $\sigma \in \Delta$ there
exists a monomial $\bm{x}^{\bm{u}} \in S_{\bm{p}}$ with $\supp(\bm{u})
\subseteq \widehat{\sigma}$.  It follows that a sufficiently large
power of each generator of $B$ belongs to the ideal $\langle
S_{\bm{p}} \rangle$ which implies $B \subseteq \sqrt{\langle
S_{\bm{p}} \rangle}$.
\end{proof}

Our motivating example of a pair $(S, B)$ is the homogeneous
coordinate ring and irrelevant ideal of a toric variety introduced in
\cite{Cox}.  Let $X$ be a simplicial toric variety over a field
$\Bbbk$ determined by a fan $\Delta$ in $\mathbb{R}^{d}$.  By
numbering the rays (one-dimensional cones) of $\Delta$, we identify
$\Delta$ with a simplicial complex on $[n]$.  We write $\bm{b}_{1},
\dotsc, \bm{b}_{n}$ for the unique minimal lattice vectors generating
the rays and we assume that the positive hull of $\scr{B} := \{
\bm{b}_{1}, \dotsc, \bm{b}_{n} \}$ is $\mathbb{R}^{d}$.  The set
$\scr{B}$ gives rise to a short exact sequence
\[
0 \longrightarrow \mathbb{Z}^{d} \xrightarrow{\;\; [ \bm{b}_{1} \;
\dotsb \; \bm{b}_{n} ]^{\textsf{T}} \;\;} \mathbb{Z}^{n}
\longrightarrow G \longrightarrow 0 \, .
\]
Tensoring this sequence with $\mathbb{R}$, we obtain
\eqref{dualconfiguration}.  Geometrically, $G$ is the divisor class
group (or Chow group) of $X$; see \S3.4 in \cite{fulton}.  The image
of the $i$th standard basis vector of $\mathbb{Z}^{n}$ in $G$ is
denoted by $\bm{a}_{i}$.  Observe that the set $\scr{A} = \{
\bm{a}_{1}, \dotsc, \bm{a}_{n} \}$ is uniquely determined up to an
automorphism of $G$.  The \emph{homogeneous coordinate ring} of $X$ is
the polynomial ring $S = \Bbbk[x_{1}, \ldots, x_{n}]$ with the
$G$\nobreakdash-grading induced by $\deg(x_{i}) = \bm{a}_{i}$ and the
\emph{irrelevant ideal} is $B = \bigl\langle \prod_{i \not\in \sigma}
x_{i} : \sigma \in \Delta \bigr\rangle$.

This geometric choice of a pair $(S, B)$ fits into the algebraic
framework developed at the start of this section if the fan of $X$
corresponds to a regular triangulation of $\scr{B}$.  By Theorem~2.6
in \cite{HauselSturmfels}, this is equivalent to saying that $X$ is
semi-projective.  In particular, this holds whenever $X$ is
projective.  Conversely, an algebraic pair $(S, B)$ is the homogeneous
coordinate ring and irrelevant ideal of a simplicial toric variety if
and only if, for each $i$, the ray $\pos(\overline{\bm{b}}_{i})$
belongs to the triangulation associated to the chamber $\Gamma$.

\begin{example} 
The pair $(S, B)$ described in Example~\ref{e:standard} corresponds to
the weighted projective space $X = \mathbb{P}(\bm{a}_{1}, \dotsc,
\bm{a}_{n})$.  In particular, when $X = \mathbb{P}^{d}$, the
homogeneous coordinate ring $S$ is the standard graded polynomial ring
and the irrelevant ideal $B$ is the unique graded maximal ideal.
\end{example}

\begin{example} \label{e:hirzvariety}
If $(S,B)$ is the pair described in Example~\ref{e:hirzebruch} where
$\Gamma = \mathbb{R}_{\geq 0}^{2}$ then the associated toric variety
is the Hirzebruch surface (or rational scroll) $\mathbb{F}_{t} =
\mathbb{P} \big( \mathscr{O}_{\mathbb{P}^{1}} \oplus
\mathscr{O}_{\mathbb{P}^{1}}(t) \big)$.
\end{example}

\begin{example} 
If $X$ is the product of projective space $\mathbb{P}^{d} \times
\mathbb{P}^{e}$, then the homogeneous coordinate ring $S =
\Bbbk[x_{0}, \dotsc, x_{d}, y_{0}, \dotsc, y_{e}]$ has the
$\mathbb{Z}^{2}$-grading induced by $\deg(x_{i}) = \left[
\begin{smallmatrix} 1 \\ 0 \end{smallmatrix} 
\right]$ and $\deg(y_{i}) = \left[
\begin{smallmatrix} 0 \\ 1 \end{smallmatrix} 
\right]$ and the irrelevant ideal $B$ is $\langle x_{0}, \dotsc, x_{d}
\rangle \cap \langle y_{0} \dotsc, y_{e} \rangle$.
\end{example}

When the pair $(S, B)$ corresponds to a simplicial toric variety $X$,
the chamber $\Gamma$ and the semigroups $\scr{K}$ and $\scr{K}^{\sat}$
have geometric interpretations.  Assuming all the maximal cones are
$d$-dimensional, we have $\Pic(X) \otimes \mathbb{Q} \cong G \otimes
\mathbb{Q}$.  As \S3 in \cite{coxrecent} indicates, the chamber
$\Gamma$ is the closure of the K\"{a}hler cone of $X$.  The dual of
the K\"{a}hler cone is the Mori cone of effective $1$-cycles modulo
numerical equivalence.  The semigroup $\scr{K}^{\sat}$ corresponds to
the numerically effective Weil divisors on $X$ up to rational
equivalence and elements in $\scr{K}$ are nef line bundles.  When $X$
is smooth, the class group $G$ is torsion-free and $\scr{K} =
\scr{K}^{\sat}$.

The next example demonstrates that $\scr{K}$ can be complicated even
when $X$ is smooth.

\begin{example}  \label{e:complicatedC}
Consider the following smooth resolution $X$ of weighted projective
space $\mathbb{P}(2,3,7,1)$.  Specifically, the set $\scr{B}$
corresponds to the columns of the matrix
\[
\left[
\begin{array}{ccccccccccc}
1 & 0 & 0 & -2 & 0 & 0 & -1 & 0 & -1 & -1 & -1 \\
0 & 1 & 0 & -3 & 0 & -1 & -1 & -1 & -1 & -2 & -2 \\
0 & 0 & 1 & -7 & -1 & -2 & -3 & -3 & -4 & -4 & -5
\end{array}
\right] \, .
\]
\begin{figure}[ht]
\psfrag{1}{\scriptsize $1$}
\psfrag{2}{\scriptsize $2 = \infty$}
\psfrag{3}{\scriptsize $3$}
\psfrag{4}{\scriptsize $4$}
\psfrag{5}{\scriptsize $5$}
\psfrag{6}{\scriptsize $6$}
\psfrag{7}{\scriptsize $7$}
\psfrag{8}{\scriptsize $8$}
\psfrag{9}{\scriptsize $9$}
\psfrag{10}{\scriptsize $10$}
\psfrag{11}{\scriptsize $11$}
\epsfig{file=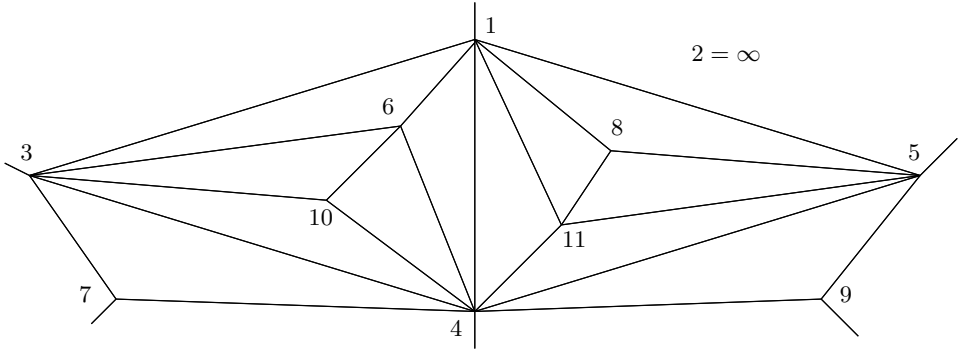, width=5in}
\caption{The triangulation for Example~\ref{e:complicatedC}
\label{f:complicatedC} }
\end{figure} 
Figure~\ref{f:complicatedC} illustrates the regular triangulation and
the irrelevant ideal is
\[
\text{\footnotesize$
\begin{array}{rllll} 
B = \langle \mspace{-15mu} &x_{1}x_{2}x_{3}x_{4}x_{6}x_{7}x_{9}x_{10}, 
&x_{2}x_{3}x_{4}x_{5}x_{6}x_{7}x_{9}x_{10}, 
&x_{1}x_{2}x_{3}x_{6}x_{7}x_{8}x_{9}x_{10}, \\
&x_{2}x_{3}x_{5}x_{6}x_{7}x_{8}x_{9}x_{10}, 
&x_{1}x_{2}x_{3}x_{5}x_{7}x_{8}x_{9}x_{11}, 
&x_{1}x_{2}x_{4}x_{5}x_{7}x_{8}x_{9}x_{11}, \\
&x_{1}x_{2}x_{5}x_{6}x_{7}x_{8}x_{9}x_{11}, 
&x_{1}x_{2}x_{3}x_{6}x_{7}x_{8}x_{10}x_{11}, 
&x_{1}x_{3}x_{4}x_{6}x_{7}x_{8}x_{10}x_{11}, \\
&x_{1}x_{3}x_{5}x_{6}x_{7}x_{8}x_{10}x_{11}, 
&x_{2}x_{3}x_{4}x_{6}x_{7}x_{9}x_{10}x_{11}, 
&x_{1}x_{2}x_{5}x_{6}x_{8}x_{9}x_{10}x_{11}, \\
&x_{1}x_{3}x_{5}x_{6}x_{8}x_{9}x_{10}x_{11}, 
&x_{1}x_{4}x_{5}x_{6}x_{8}x_{9}x_{10}x_{11}, 
&x_{2}x_{3}x_{5}x_{7}x_{8}x_{9}x_{10}x_{11}, \\
&x_{2}x_{4}x_{5}x_{7}x_{8}x_{9}x_{10}x_{11}, 
&x_{3}x_{4}x_{6}x_{7}x_{8}x_{9}x_{10}x_{11}, 
&x_{4}x_{5}x_{6}x_{7}x_{8}x_{9}x_{10}x_{11} \rangle \, .
\end{array}$}
\]
In this case, $\Gamma$ is an $8$-dimensional cone with $16$ extremal
rays.  Using \cite{normaliz}, we determine that set of minimal
generators for $\scr{K}$ has 25 elements all of which lie on the
boundary of $\Gamma$.
\end{example}

\section{A Topological Formula for Local Cohomology}

Throughout this paper, $M$ denotes a finitely generated $G$-graded
$S$-module.  In this section we derive a vanishing theorem for the
local cohomology modules $H_{B}^{i}(M)$ from a topological formula for
$H_{B}^{i}(S)$.  We refer to \cite{BrodmannSharp} for background
information on local cohomology.  A module $M$ is \emph{$B$-torsion}
if $M = H_{B}^{0}(M)$.  If $M$ is $B$-torsion then $H_{B}^{i}(M) = 0$
for $i > 0$.  For an element $g \in S$, we set
\[
( 0 :_{M} g) := \{ f \in M : gf = 0 \} = \Ker \bigl( M
\xrightarrow{\;\; g \;\;} M \bigl) \, .
\]  
This submodule is zero when $g$ is a nonzerodivisor on $M$.  We say an
element $g \in S$ is \emph{almost a nonzerodivisor on $M$} if $(0:_{M}
g)$ is a $B$-torsion module.

\begin{proposition} \label{p:nonzerodivisor}
Let $\Bbbk$ be an infinite field.  If $\bm{p} \in \scr{K}$ and $g \in
S$ is a sufficiently general form of degree $\bm{p}$, then $g$ is
almost a nonzerodivisor on $M$.  In other words, there is an open
dense set of degree $\bm{p}$ forms $g$ for which $(0 :_{M} g)$ is
$B$-torsion.
\end{proposition}

\begin{proof}
The module $M' := (0:_{M} g)$ is $B$-torsion if each element in $M'$
is annihilated by some power of $B$.  This is equivalent to saying
that the localization $M_{P}' = 0$ for all prime ideals $P$ in $S$ not
containing $B$.  In other words, $g$ is a nonzerodivisor on the
localization $M_{P}$.  Hence, it suffices to show that $g$ is not
contained in any of the associated primes of $M$ except possibly those
which contain $B$.

To accomplish this, observe that each prime ideal $P$ in $S$ which
does not contain $B$ intersects $S_{\bm{p}}$ in a proper subspace.
Otherwise $\langle S_{\bm{p}} \rangle \subseteq P$ and
Lemma~\ref{l:BandC} implies that $B \subseteq \sqrt{\langle S_{\bm{p}}
\rangle} \subseteq \sqrt{P} = P$ which contradicts the hypothesis on
$P$.  Because $M$ has only finitely many associated primes, our
observation shows that $g \in S_{\bm{p}}$ is almost a nonzerodivisor
on $M$ if it lies outside a certain finite union of proper subspaces.
Since $\Bbbk$ is infinite, the vector space $S_{\bm{p}}$ is not a
finite union of proper subspaces.
\end{proof}

To give a uniform vanishing result for local cohomology, we assume for
the remainder of this section that $\pos(\overline{\scr{A}})$ is a
pointed cone with $\overline{\bm{a}}_{i} \neq \bm{0}$ for each $i$ or
equivalently that $\overline{\scr{A}}$ is an acyclic vector
configuration; see \S6.2 in \cite{Ziegler}.  This condition holds if
$S$ is the homogeneous coordinate ring of a projective toric variety.
We can rephrase this assumption by saying that the dual configuration
$\overline{\scr{B}}$ is totally cyclic or $\pos(\overline{\scr{B}}) =
\mathbb{R}^{d}$.  Hence, any regular triangulation of
$\overline{\scr{B}}$ is a complete fan in $\mathbb{R}^{d}$.  We may
identify $\Delta$ with a simplicial $(d-1)$-sphere and the collection
$\widehat{\Delta} := \Delta \cup \{ [n] \}$ (ordered by inclusion)
with the face poset of a finite regular cell decomposition of a
$d$-ball.  It follows that there is an incidence function
$\varepsilon$ on $\widehat{\Delta}$; see \S6.2 in \cite{BH}.

For $\sigma \subseteq [n]$, let $\bm{x}_{\sigma}$ denote the
squarefree monomial $\prod_{i \in \sigma} x_{i} \in S$.  The
\emph{canonical \v{C}ech complex} associated to $\widehat{\Delta}$ is
the following chain complex:
\[
\mathbf{\check{C}} := 0 \longrightarrow \check{C}^{0} \xrightarrow{\;\;
\partial^{0} \;\;} \check{C}^{1} \xrightarrow{\;\; \partial^{1} \;\;}
\check{C}^{2} \xrightarrow{\;\; \partial^{2} \;\;} \dotsb
\xrightarrow{\;\; \partial_{d-1} \;\;} \check{C}^{d} \xrightarrow{\;\;
\partial_{d} \;\;} \check{C}^{d+1} \longrightarrow 0 \, ,
\]
where $\check{C}^{0} := S$, $\check{C}^{i} := \bigoplus_{\sigma \in
\Delta, |\sigma| = d+1-i} S[\bm{x}_{\widehat{\sigma}}^{-1}]$ for $i >
0$ and $\partial^{i} \colon \check{C}^{i} \longrightarrow
\check{C}^{i+1}$ is composed of homomorphisms $\varepsilon(\sigma,
\tau) \cdot \text{nat} \colon S[\bm{x}_{\widehat{\sigma}}^{-1}]
\longrightarrow S[\bm{x}_{\widehat{\tau}}^{-1}]$.  By combining
Corollary~2.13 in \cite{BS2} and Theorem~6.2 in \cite{Miller}, we
deduce that $H_{B}^{i}(M) = H^{i}(\mathbf{\check{C}} \otimes_{S} M)$
for all $S$-modules $M$.  It follows that $H_{B}^{i}(M) = 0$ for all
$i < 0$ or $i > d + 1$.

The next proposition gives a combinatorial description for the local
cohomology of the polynomial ring $S$.  We will use this formula to
compute the multigraded regularity in some examples.  This result is a
variant on the formulae found in \cite{EMS} and \cite{RWY}.  For
$\sigma \in \Delta$, let $\Delta_{\sigma}$ be the induced subcomplex
$\{ \tau \in \Delta : \tau \subseteq \sigma \}$ of $\Delta$.  The
modules $H_{B}^{i}(S)$ have a $\mathbb{Z}^{n}$-grading which refines
their $G$-grading.

\begin{proposition} \label{p:firsttopoprop}
Let $\widetilde{H}^{i}(-)$ denote the $i$th reduced cellular
cohomology group with coefficients in $\Bbbk$.  If $\bm{u} \in
\mathbb{Z}^{n}$ and $\sigma := \negg(\bm{u})=\{ j : \bm{u}_j <0 \}$,
then we have
\[ 
H_{B}^{i}(S)_{\bm{u}} \cong
\begin{cases}
\widetilde{H}^{i-2} \bigl( \Delta_{\sigma} \bigr) & \text{for $i \neq
1;$} \\ 0 & \text{for $i = 1$.}
\end{cases}
\]
\end{proposition}

\begin{proof}
One easily checks (see Lemma~5.3.6 in \cite{BH}) that $\dim_{\Bbbk} S[
\bm{x}_{\widehat{\tau}}^{-1}]_{\bm{u}} = 1$ when $\tau \subseteq
\widehat{\sigma}$ and $\dim_{\Bbbk} S[
\bm{x}_{\widehat{\tau}}^{-1}]_{\bm{u}} = 0$ when $\tau \not\subseteq
\widehat{\sigma}$.  It follows that the $\bm{u}$th graded component of
$\mathbf{\check{C}}$ is isomorphic to a shift of the augmented
oriented chain complex of $\widehat{\Delta}_{\widehat{\sigma}}$ (see
\S6.2 in \cite{BH}).  More precisely, we have
$H^{i}(\mathbf{\check{C}}_{\bm{u}}) = \widetilde{H}_{d-i} (
\widehat{\Delta}_{\widehat{\sigma}} )$.  To complete the proof, we
analyze three cases:
\begin{itemize}
\item If $\sigma$ is a proper nonempty subset of $[n]$, then the
subcomplex $\widehat{\Delta}_{\widehat{\sigma}}$ equals
$\Delta_{\widehat{\sigma}}$.  Alexander duality (Theorem~71.1 in
\cite{munkres}) shows that $\widetilde{H}_{d-i} (
\Delta_{\widehat{\sigma}} )$ is isomorphic to
$\widetilde{H}^{i-2}(\Delta \setminus \Delta_{\widehat{\sigma}})$.
Because $\Delta_{\sigma}$ is a deformation retract of $\Delta
\setminus \Delta_{\widehat{\sigma}}$ (Lemma~70.1 in \cite{munkres}),
we have $\widetilde{H}_{d-i} ( \Delta_{\widehat{\sigma}} ) \cong
\widetilde{H}^{i-2} ( \Delta_{\sigma} )$.  Note that
$\widetilde{H}^{i-2} ( \Delta_{\sigma} ) = 0$ when $i = 1$.
\item If $\sigma = \emptyset$, then $\widehat{\sigma} = [n]$ and
$\widehat{\Delta}_{{\widehat{\sigma}}}$ is a $d$-ball which implies
that $\widetilde{H}_{d-i} ( \widehat{\Delta}_{\widehat{\sigma}} ) = 0$
for all $i$.  It follows that $\widetilde{H}_{d-i} (
\widehat{\Delta}_{[n]} )$ is isomorphic to
$\widetilde{H}^{i-2}(\Delta_{\emptyset})$ for $i \neq 1$ and equals
$0$ for $i = 1$.
\item If $\sigma = [n]$, then $\widehat{\sigma} = \emptyset$ and
$\widehat{\Delta}_{{\widehat{\sigma}}} = \emptyset$ which implies that
$\widetilde{H}_{d-i} ( \widehat{\Delta}_{\widehat{\sigma}} ) = 0$ for
$i \neq d+1$ and $\widetilde{H}_{-1} (
\widehat{\Delta}_{\widehat{\sigma}} ) \cong \Bbbk$.  Since
$\Delta_{[n]} = \Delta$ is a $(d-1)$-sphere, we also have
$\widetilde{H}_{d-i} ( \widehat{\Delta}_{\emptyset} ) \cong
\widetilde{H}^{i-2} ( \Delta_{[n]} )$ for all $i$. \qedhere
\end{itemize}
\end{proof}

\begin{remark} \label{r:whenvanish}
Proposition~\ref{p:firsttopoprop} implies that $H_{B}^{i}(S) \neq 0$
only if $2 \leq i \leq d+1$.  Since $H_{B}^{d+1}(S)_{\bm{u}}=\Bbbk$
when $\negg(\bm{u}) = [n]$, we have $H_{B}^{d+1}(S)\neq 0$.  Moreover,
$H_{B}^{d+1}(S)$ is the only nonvanishing local cohomology module if
and only if every proper subcomplex of $\Delta$ is contractible.  This
happens precisely when $\Delta$ is the boundary of the standard
simplex.  Hence, $H_{B}^{i}(S) = 0$ for all $i \neq d+1$ if and only
if $B = \langle x_{1}, \dotsc, x_{n} \rangle$.
\end{remark}

Using Proposition~\ref{p:firsttopoprop}, we can describe the Hilbert
series of the modules $H_{B}^{i}(S)$.

\begin{corollary} \label{c:hilbertseries}
For all $i \neq 1$, we have  
\begin{equation} \label{hilbertseries:1}
\sum_{\bm{p} \in G} \dim_{\Bbbk} H_{B}^{i}(S)_{\bm{p}} \cdot
\bm{t}^{\bm{p}} = \sum_{\sigma \subseteq [n]} \frac{\dim_{\Bbbk}
\widetilde{H}^{i-2}(\Delta_{\sigma}) \cdot \prod_{j \in \sigma}
\bm{t}^{-\bm{a}_{j}}}{\prod_{j \in \widehat{\sigma}} (1-
\bm{t}^{\bm{a}_{j}}) \prod_{j \in \sigma} (1- \bm{t}^{\bm{-a}_{j}})}
\, .
\end{equation}
\end{corollary}

\begin{proof}
For $i \neq 1$, Proposition~\ref{p:firsttopoprop} implies that
\begin{multline*}
\sum_{\bm{u} \in \mathbb{Z}^{n}} \dim_{\Bbbk} H_{B}^{i}(S)_{\bm{u}}
\cdot \bm{x}^{\bm{u}}
= \sum_{\bm{u} \in \mathbb{Z}^{n}} \dim_{\Bbbk} \widetilde{H}^{i-2} 
\bigl( \Delta_{\negg(\bm{u})} \bigr) \cdot \bm{x}^{\bm{u}} \\
= \sum_{\sigma \subseteq [n]} \left( \dim_{\Bbbk}
\widetilde{H}^{i-2}(\Delta_{\sigma}) \sum_{\negg(\bm{u}) = \sigma}
\bm{x}^{\bm{u}} \right) 
= \sum_{\sigma \subseteq [n]} \frac{\dim_{\Bbbk}
\widetilde{H}^{i-2}(\Delta_{\sigma}) \cdot \prod_{j \in \sigma}
x_{j}^{-1}}{\prod_{j \in \widehat{\sigma}} (1- x_{j}) \prod_{j \in
\sigma} (1- x_{j}^{-1})} \, .
\end{multline*}
Mapping $x_{j}$ to $\bm{t}^{\bm{a}_{j}}$ establishes the corollary.
\end{proof}

We illustrate this corollary with the following example.

\begin{example} \label{e:hirz2}
If $S$ is the homogeneous coordinate ring of the Hirzebruch surface
$\mathbb{F}_{t}$ as in Examples~\ref{e:hirzebruch} \&
\ref{e:hirzvariety}, then the only non-contractible subcomplexes are
$\Delta_{\emptyset}$, $\Delta_{\{1, 3\}}$, $\Delta_{\{2, 4\}}$ and
$\Delta_{[4]}$.  Since $\widetilde{H}^{i}(\Delta_{\emptyset}) = 0$ for
all $i \neq -1$ and $\widetilde{H}^{-1}(\Delta_{\emptyset}) = \Bbbk$,
this subcomplex does not contribute to any local cohomology module.
Because $\widetilde{H}^{0}(\Delta_{\{1, 3\}}) = \Bbbk =
\widetilde{H}^{0}(\Delta_{\{2, 4\}})$, the degrees $\bm{p} \in
\mathbb{Z}^{2}$ for which $H_{B}^{2}(S)_{\bm{p}} \neq 0$ correspond to
the lattices points in the two cones $- \bm{a}_{1} - \bm{a}_{3} +
\pos(-\bm{a}_{1}, \bm{a}_{2}, - \bm{a}_{3}, \bm{a}_{4})$ and $-
\bm{a}_{2} - \bm{a}_{4} + \pos(\bm{a}_{1}, -\bm{a}_{2}, \bm{a}_{3}, -
\bm{a}_{4})$.  We also have $\widetilde{H}^{1}(\Delta_{[n]}) = \Bbbk$,
so the degrees $\bm{p} \in \mathbb{Z}^{2}$ for which
$H_{B}^{2}(S)_{\bm{p}} \neq 0$ correspond to the lattices points in
the cone $- \bm{a}_{1} - \bm{a}_{2} - \bm{a}_{3} - \bm{a}_{4} +
\pos(-\bm{a}_{1}, -\bm{a}_{2}, -\bm{a}_{3}, - \bm{a}_{4})$.  These
cones are indicated by the shaded areas in
Figure~\ref{f:nonvanishingcones} when $t = 2$.
\begin{figure}[ht]
\psfrag{H1}{$H^2_B(S)$}
\psfrag{H2}{$H^3_B(S)$}
\epsfig{file=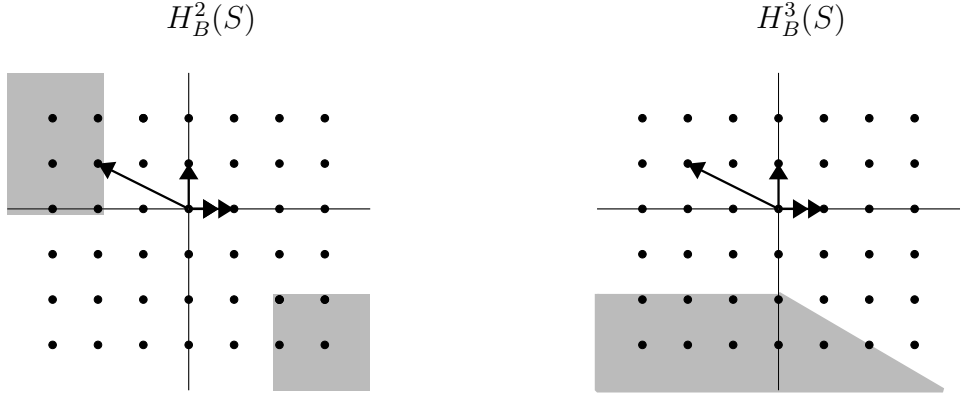, width=5in}
\caption{Degrees $\bm{p} \in \mathbb{Z}^{2}$ for which
$H_{B}^{i}(S)_{\bm{p}} \neq 0$ in Example
\ref{e:hirz2}. \label{f:nonvanishingcones}}
\end{figure}
\end{example}

The next result extends a well-known vanishing theorem
for ample line bundles on a complete toric variety; see
\eqref{localglobal} for the explicit connection.  

\begin{corollary} \label{t:vanishing} 
If $\bm{p}$ belongs to $\scr{K}^{\sat}$, then $H_{B}^{i}(S)_{\bm{p}}$
vanishes.  
\end{corollary}

\begin{proof}
Fix $\bm{p} \in \scr{K}^{\sat}$ and suppose $\bm{u} \in
\mathbb{Z}^{n}$ satisfies $\deg(\bm{x}^{\bm{u}}) = \bm{p}$.  The
vector $\bm{u}$ defines a function on $\psi_{\bm{u}} \colon
\mathbb{R}^{d} \longrightarrow \mathbb{R}$ that is linear on the cone
$\pos(\overline{\bm{b}}_{i} : i \in \sigma)$ for all $\sigma \in
\Delta$ and satisfies $\psi_{\bm{u}}(\overline{\bm{b}}_{i}) = - u_{i}$
for $1 \leq i \leq n$.  Lemma~3.2 in \cite{BGS} implies that $\bm{p}
\in \scr{K}^{\sat}$ if and only if $\psi_{\bm{u}}$ is convex which
means $\psi_{\bm{u}}(\bm{w} + \bm{w}') \geq \psi_{\bm{u}}(\bm{w}) +
\psi_{\bm{u}}(\bm{w}')$ for all $\bm{w}$, $\bm{w}' \in
\mathbb{R}^{d}$.  Hence, the set $\{ \bm{w} \in \mathbb{R}^{d} :
\psi_{\bm{u}}(\bm{w}) \geq 0 \}$ is convex.  It follows that
simplicial complex $\Delta_{\negg(\bm{u})}$ corresponds to the
intersection of a convex set with the $(d-1)$-sphere.  Therefore,
$\Delta_{\negg(\bm{u})}$ is contractible and
Proposition~\ref{p:firsttopoprop} implies that $H_{B}^{i}(S)_{\bm{p}}
= 0$ which proves the statement.
%
\end{proof}
%

As a further corollary, we obtain a vanishing theorem for the local
cohomology of any finitely generated $G$-graded $S$-module.
Geometrically, this result corresponds to Fujita's vanishing theorem
(Theorem~1.4.32 in \cite{lazarsfeld}).  We first record a useful
observation.

\begin{remark} \label{r:nonempty}
Let $\scr{D}$ be a subsemigroup of $G$ consisting of all points
$\bm{p} \in G$ such that their image $\overline{\bm{p}}$ in
$\mathbb{R}^{r}$ lies in a fixed convex cone $C$.  If $\bm{p} +
\scr{D}$ and $\bm{q} + \scr{D}$ are two shifts of $\scr{D}$, then
their intersection is nonempty.  In fact, if $\bm{p} \in \scr{D}$,
with $\overline{\bm{p}} \in \interior C$, then for all sufficiently
large $j \in \mathbb{N}$ we have $j \bm{p} \in \bm{q} + \scr{D}$ for
any $\bm{q} \in G$.
\end{remark}

\begin{corollary} \label{c:fujita}
There is $\bm{m} \in \scr{K}^{\sat}$ such that $H_{B}^{i}(M)_{\bm{p}}
= 0$ for all $i$ and all $\bm{p} \in \bm{m} + \scr{K}^{\sat}$.  In
particular, $H_{B}^{i}(M)_{\bm{p}}$ vanishes for all $i$ and all
$\bm{p} \in \scr{K}^{\sat}$ sufficiently far into the interior of
$\scr{K}^{\sat}$.  If desired we may assume $\bm{m} \in \scr{K}$.
\end{corollary}

\begin{proof}
We proceed by induction on the projective dimension $\pdim(M)$ of $M$.
Since $M$ is finitely generated, there is a short exact sequence $0
\longrightarrow E_{1} \longrightarrow E_{0} \longrightarrow M
\longrightarrow 0$ where $E_{0} = \bigoplus_{1 \leq j \leq h}
S(-\bm{q}_{j})$ for some $\bm{q}_{j} \in G$.  The associated long
exact sequence contains
\begin{equation} \label{fujita1}
H_{B}^{i}(E_{0})_{\bm{p}} \longrightarrow H_{B}^{i}(M)_{\bm{p}}
\longrightarrow H_{B}^{i+1}(E_{1})_{\bm{p}} \, .
\end{equation}
If we choose $\bm{m}' \in \bigcap_{1 \leq j \leq h} \bigl( \bm{q}_{j}
+ \scr{K}^{\sat} \bigr)$, which is possible by
Remark~\ref{r:nonempty}, then Theorem~\ref{t:vanishing} implies that
the left module in \eqref{fujita1} vanishes for all $i$ and all
$\bm{p} \in \bm{m}' + \scr{K}^{\sat}$.  When $\pdim(M) = 0$, we have
$E_{1} = 0$ and there is nothing more to prove.  Otherwise, we have
$E_{1} \neq 0$ and $\pdim(E_{1}) < \pdim(M)$.  Hence, the induction
hypothesis provides $\bm{m}'' \in G$ such that the right module in
\eqref{fujita1} vanishes for all $\bm{p} \in \bm{m}'' +
\scr{K}^{\sat}$.  Therefore, by choosing $\bm{m} \in \bigl( \bm{m}' +
\scr{K}^{\sat} \bigr) \cap \bigl( \bm{m}'' + \scr{K}^{\sat} \bigr)$,
which is again possible by Remark~\ref{r:nonempty}, we see that the
middle module in \eqref{fujita1} vanishes for all $i$ and all $\bm{p}
\in \bm{m} + \scr{K}^{\sat}$.

Lastly, because $\scr{K}^{\sat}$ is the saturation of $\scr{K}$, there
exists some $\bm{m}''' \in \scr{K} \cap (\bm{m}+\scr{K}^{\sat})$.
Since $\bm{m}''' + \scr{K}^{\sat} \subseteq \bm{m}+\scr{K}^{\sat}$,
the corresponding vanishing statement holds for $\bm{m}'''$ which
establishes the last part of the corollary.
\end{proof}

We end this section by examining the assumption that
$\pos(\overline{\scr{A}})$ is acyclic.  Firstly, we can remove the
condition that $\bm{a}_{i} \neq \bm{0}$ for all $i$.

\begin{remark}
Let $\sigma := \{ i \in [n] : \overline{\bm{a}}_{i} = \bm{0} \}$.  For
each $j \in \sigma$ the variable $x_{j}$ does not divide any minimal
generator of $B$.  In other words, the set $\sigma$ is contained in
every facet of $\Delta$.  Hence, the simplicial complex $\Delta$ is a
cone over the induced subcomplex $\Delta_{\widehat{\sigma}}$.  By
replacing $\Delta$ with $\Delta_{\widehat{\sigma}}$, we can extend
Proposition~\ref{p:firsttopoprop} and its corollaries to this more
general situation.
\end{remark}

On the other hand, the next example shows that we cannot eliminate the
hypothesis that $\pos(\overline{\scr{A}})$ is a pointed cone.

\begin{example} \label{e:novanishing}
If $G = \mathbb{Z}$ and $\scr{A} = \{ \bm{1}, \bm{-1}, \bm{1} \}$,
then $S = \Bbbk[x_{1}, x_{2}, x_{3}]$ has the
$\mathbb{Z}$\nobreakdash-grading induced by $\deg(x_{2}) = \bm{-1}$
and $\deg(x_{1}) = \deg(x_{3}) = \bm{1}$.  Choosing the chamber
$\Gamma$ to be $\mathbb{R}_{\geq 0}$, we have $\scr{K} =
\scr{K^{\sat}} = \mathbb{N}$.  The corresponding monomial ideal is $B
= \langle x_{1}, x_{3} \rangle$ and the dual vector configuration
corresponds to the columns of the matrix $\left[
\begin{smallmatrix} 1 & 1 & 0 \\ 0 & 1 & 1 \end{smallmatrix} 
\right]$.  Since the first paragraph of the proof of
Proposition~\ref{p:firsttopoprop} applies in this situation,
$H_{B}^{2}(S)_{\bm{u}} = \widetilde{H}_{0} (
\widehat{\Delta}_{\widehat{\sigma}} )$ where $\sigma = \negg(\bm{u})$.
If $\sigma = \{ 2 \}$, then $\widehat{\sigma} = \{ 1,3 \}$ and
$\widehat{\Delta}_{\{1,3\} } = \bigl\{ \emptyset, \{ 1 \}, \{ 3 \}
\bigr\}$.  It follows that $\widetilde{H}_{0} (
\widehat{\Delta}_{\widehat{\sigma}} ) \neq 0$ which implies that
$H^2_B(S)_{\bm{u}} \neq 0$ when $\negg(\bm{u}) = \{2 \}$.  However,
for all $\bm{p} \in \mathbb{N}$, there exists a $\bm{u} \in \mathbb
Z^3$ with $\negg(\bm{u})= \{2 \}$.  We conclude that
$H_{B}^{2}(S)_{\bm{p}} \neq 0$ for all $\bm{p} \in \scr{K}$.
\end{example}

\section{Regularity of $S$-modules}

In this section, we define the multigraded regularity of a $G$-graded
$S$-module.  As in the standard case, multigraded regularity involves
the vanishing of certain graded components of local cohomology
modules.  We show that, in certain cases, the definition is equivalent
to an apparently weaker vanishing statement.

Before giving the definition of multigraded regularity, we collect
some notation.  Throughout this paper, $\scr{C} := \{ \bm{c}_{1},
\dotsc, \bm{c}_{\ell} \}$ will be a fixed subset of $G$.  We write
$\mathbb{N} \scr{C}$ for the subsemigroup of $G$ generated by
$\scr{C}$.  A subset $\scr{D}$ of $G$ is called an \emph{$\mathbb{N}
\scr{C}$\nobreakdash-module} if for all $\bm{p} \in \scr{D}$ and all
$\bm{q} \in \mathbb{N} \scr{C}$ we have $\bm{p} + \bm{q} \in \scr{D}$.
If $\scr{D} \subseteq G$ and $i \in \mathbb{Z}$, then $\scr{D}[i] :=
\bigcup \bigl( \tfrac{i}{|i|} ( \lambda_{1} \bm{c}_{1} + \dotsb +
\lambda_{\ell} \bm{c}_{\ell}) + \scr{D} \bigr) \subseteq G$ where the
union is over all $\lambda_{1}, \dotsc, \lambda_{\ell} \in \mathbb{N}$
such that $\lambda_{1} + \dotsb + \lambda_{\ell} = |i|$.  For $\bm{p}
\in G$, we clearly have $\bm{p} + \scr{D}[i] = (\bm{p} + \scr{D})[i]$.
Moreover, if $\scr{D}$ is an $\mathbb{N} \scr{C}$-module then
$\scr{D}[i]$ is also an $\mathbb{N} \scr{C}$-module and $\scr{D}[i+1]
\subseteq \scr{D}[i]$.  The following definition includes
Definition~\ref{i:def} as a special case.

\begin{definition} \label{d:regdefn}
If $\bm{m} \in G$, then the $S$-module $M$ is \emph{$\bm{m}$-regular}
(with respect to $\scr{C}$) if $H_{B}^{i}(M)_{\bm{p}} = 0$ for all $i
\geq 0$ and all $\bm{p} \in \bm{m} + \mathbb{N} \scr{C}[1-i]$.  The
\emph{regularity of $M$}, denoted $\reg(M)$, is the set $\{ \bm{m} \in
G : \text{$M$ is $\bm{m}$-regular} \}$.
\end{definition}

When $S$ has the standard grading, $B = \langle x_{1}, \dotsc, x_{n}
\rangle$, and $\scr{C} = \{ \bm{1} \}$, the definition is simply: $M$
is $\bm{m}$-regular if and only if $H_{B}^{i}(M)_{\bm{p}} = 0$ for all
$i \geq 0$ and all $\bm{p} \geq \bm{m}-i+1$.  This is equivalent to
the standard definition of Castelnuovo-Mumford regularity; for example
see page~282 of \cite{BrodmannSharp}.  Moreover, if $S$ has an
$\mathbb{N}$-grading as in Example~\ref{e:standard}, $B = \langle
x_{1}, \dotsc, x_{n} \rangle$ and $\scr{C} = \{ \bm{1} \}$, then the
definition of $\bm{m}$-regular is compatible with Definition~4.1 in
\cite{Benson}.

Our definition for $\reg(M)$ conflicts with the standard notation used
for Castelnuovo-Mumford regularity.  Traditionally, $G = \mathbb{Z}$
and the regularity of $M$ referred to the smallest $\bm{m} \in G$ such
that $M$ is $\bm{m}$-regular.  In the multigraded setting, the group
$G$ is not equipped with a natural ordering so one cannot choose a
smallest degree $\bm{m}$.  More importantly, the set $\reg(M)$ may not
even be determined by a single element of $G$.
Example~\ref{e:hirzreg1} illustrates this phenomenon.  For these
reasons, we regard the regularity of $M$ as a subset of $G$.  From
this vantage point, giving a bound on the smallest $\bm{m}$ such that
$M$ is $\bm{m}$-regular should be interpreted as describing a subset
of $\reg(M)$.

\begin{example} \label{e:standardreg}
Let $S = \Bbbk[x_{1}, \dotsc, x_{n}]$ have the $\mathbb{Z}$-grading
induced by $\deg(x_{i}) = \bm{a}_{i} > 0$ and the ideal $B$ be
$\langle x_{1}, \dotsc, x_{n} \rangle$ as in Example~\ref{e:standard}.
If $\scr{C} = \{ \bm{1} \}$, then Remark~\ref{r:whenvanish} implies
that $\reg(S) = \{ \bm{m} \in \mathbb{Z} : \bm{m} \geq n - \bm{a}_{1}
- \dotsb - \bm{a}_{n} \}$.  On the other hand, if $\bm{a}_{i} = 1 $
for $1 \leq i \leq n$ and $\scr{C} = \{ \bm{c}_{1}, \dotsc,
\bm{c}_{\ell} \}$ with $0 < \bm{c}_{1} \leq \dotsb \leq
\bm{c}_{\ell}$, then we have $\reg(S) = \bigl\{ \bm{m} \in \mathbb{Z}
: \bm{m} \geq (n-1) (\bm{c}_{\ell} - 1) \bigr\}$.
\end{example}

\begin{example} \label{e:hirzreg1}
Let $t \in \mathbb{N}$. Suppose that $S$ is the homogeneous coordinate
ring of the Hirzebruch surface $\mathbb{F}_{t}$; see
Example~\ref{e:hirzvariety}.  If $\scr{C} = \left\{ \left[
\begin{smallmatrix} 1 \\ 0 \end{smallmatrix} 
\right], \left[
\begin{smallmatrix} 0 \\ 1 \end{smallmatrix} 
\right] \right\}$ then Example~\ref{e:hirz2} implies that
\[
\reg(S) = 
\begin{cases}
\mathbb{N}^{2} & \text{for $t = 0$, $1$;} \\ 
\bigl( \left[ 
\begin{smallmatrix} t-1 \\ 0 \end{smallmatrix} 
\right] + \mathbb{N}^{2} \bigr) \cup \bigl( \left[
\begin{smallmatrix} 0 \\ 1 \end{smallmatrix} 
\right] + \mathbb{N}^{2} \bigr) & \text{for $t \geq 2$.}
\end{cases}
\]
Observe that when $t \geq 2$ the set $\reg(S)$ is not determined by a
single element of $G$.  Moreover, we have $\interior \scr{K}^{\sat}
\subset \reg(S)$ for all $t \in \mathbb{N}$, but $\mathbf{0} \not\in
\reg(S)$ for $t \geq 2$.  The shaded region in Figure~\ref{f:hirzreg}
represents $\reg(S)$ when $t = 2$.
\begin{figure}[ht]
\epsfig{file=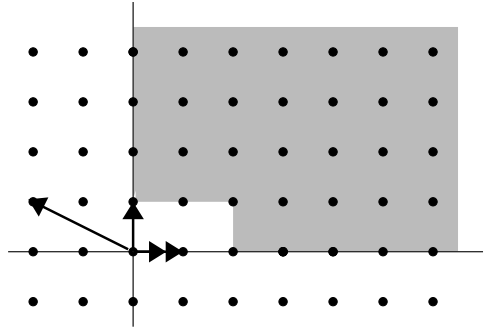, width=2.5in}
\caption{The $\reg(S)$ for the homogeneous coordinate ring $S$ of
$\mathbb{F}_{2}$. \label{f:hirzreg}}
\end{figure}
\end{example}

Without additional hypotheses, it is possible that $\reg(M)$ is the
empty set.  Indeed, if $S = \Bbbk[x_{1}, x_{2}, x_{3}]$ has the
$\mathbb{Z}$-grading described in Example~\ref{e:novanishing}, then
$\reg(S) = \emptyset$.  Fortunately, there is a large class of pairs
$(S, B)$ for which $\reg(M) \neq \emptyset$ for all finitely generated
$S$-modules $M$.  Specifically, we have the following:

\begin{proposition}
If $\pos(\overline{\scr{A}})$ is a pointed cone and $\mathbb{N}
\scr{C} \cap \interior \scr{K}^{\sat} \neq \emptyset$, then every
module $M$ is $\bm{m}$-regular for some $\bm{m} \in G$.
\end{proposition}

\begin{proof}
Corollary~\ref{c:fujita} states that there is $\bm{m} \in
\scr{K}^{\sat}$ such that $H_{B}^{i}(M)_{\bm{p}} = 0$ for all $i$ and
all $\bm{p} \in \bm{m} + \scr{K}^{\sat}$.  By hypothesis, there exists
$\bm{c} \in \mathbb{N} \scr{C} \cap \interior \scr{K}^{\sat}$.
Remark~\ref{r:nonempty} implies that there exists $k \in \mathbb{N}$
such that
\[ 
k \bm{c} \in \bigcap\limits_{\begin{subarray}{c} 
\lambda_{1}, \dotsc, \lambda_{\ell} \in \mathbb{N} \\
\lambda_{1} + \dotsb + \lambda_{\ell} \leq d
\end{subarray}} 
( \bm{m} + \lambda_{1} \bm{c}_{1} + \dotsb + \lambda_{\ell}
\bm{c}_{\ell} + \scr{K}^{\sat} ) \, .
\]
We conclude that $k \bm{c} \in \reg(M)$.
\end{proof}

To show that a module $M$ is $\bm{m}$-regular directly from
Definition~\ref{d:regdefn}, one must verify that an infinite number of
graded components of the local cohomology vanish.  In fact, under the
appropriate hypothesis, one need only check a finite number of
components.  To prove this result, we introduce the following two
weaker versions of regularity.

\begin{definition} \label{d:weakreg}
Given $k \in \mathbb{N}$, the module $M$ is \emph{$\bm{m}$-regular
from level $k$} if $H_{B}^{i}(M)_{\bm{p}} = 0$ for all $i \geq k$ and
all $\bm{p} \in \bm{m} + \mathbb{N} \scr{C}[1-i]$.  In particular, $M$
is $\bm{m}$-regular if and only if it is $\bm{m}$-regular from level
$0$.  We set $\reg^{k}(M) := \{ \bm{m} \in G : \text{$M$ is
$\bm{m}$-regular from level $k$} \}$.  For $k > 0$, we say $M$ is
\emph{weakly} $\bm{m}$-regular from level $k$ if
$H_{B}^{i}(M)_{\bm{p}} = 0$ for all $i \geq k$ and all $\bm{p} =
\bm{m} - \lambda_{1} \bm{c}_{1} - \dotsb - \lambda_{\ell}
\bm{c}_{\ell}$ where $\lambda_{j} \in \mathbb{N}$ and $\lambda_{1} +
\dotsb + \lambda_{\ell} = i - 1$.
\end{definition}

Our goal is to show that, when $\scr{C} \subseteq \scr{K}$, $M$ is
$\bm{m}$-regular if and only if it is weakly $\bm{m}$-regular from
level $1$ and $H_{B}^{0}(M)_{\bm{p}} = 0$ for $\bm{p} \in \bigcup_{1
\leq j \leq \ell} (\bm{m} + \bm{c}_{j} + \mathbb{N} \scr{C})$.  To
accomplish this, we need the following fact.

\begin{lemma} \label{l:weakly}
Let $k$ be a positive integer.  If $M$ is weakly $\bm{m}$-regular from
level $k$ and $g \in S_{\bm{c}_{j}}$ is almost a nonzerodivisor on $M$
for some $1 \leq j \leq \ell$, then $M/g_{j}M$ is also weakly
$\bm{m}$-regular from level $k$.
\end{lemma}

\begin{proof}
Since $g_{j}$ is almost a nonzerodivisor on $M$, the submodule $(0
:_{M} g_{j})$ is a $B$-torsion module.  We set $\overline{M} := M/ (0
:_{M} g_{j})$.  Because $H_{B}^{i}\bigl( (0 :_{M} g_{j}) \bigr) = 0$
for $i > 0$, the long exact sequence associated to the exact sequence
$0 \longrightarrow (0 :_{M} g_{j}) \longrightarrow M \longrightarrow
\overline{M} \longrightarrow 0$ implies that $H_{B}^{i}(M) =
H_{B}^{i}(\overline{M})$ for $i > 0$.  Since $M$ is weakly
$\bm{m}$-regular from level $k$ and $k > 0$, it follows that
$\overline{M}$ is weakly $\bm{m}$-regular from level $k$.  On the
other hand, the long exact sequence associated to the exact sequence
$0 \longrightarrow \overline{M}(- \bm{c}_{j}) \longrightarrow M
\longrightarrow M/g_{j}M \longrightarrow 0$ contains
\begin{equation} \label{weakly:1}
H_{B}^{i}(M)_{\bm{p}} \longrightarrow
H_{B}^{i}(M/g_{j}M)_{\bm{p}} \longrightarrow
H_{B}^{i+1}(\overline{M})_{\bm{p}-\bm{c}_{j}} \, .
\end{equation}  
Since $M$ and $\overline{M}$ are weakly $\bm{m}$-regular from level
$k$, the left and right modules in \eqref{weakly:1} vanish when $i
\geq k$ and $\bm{p} = \bm{m} - \lambda_{1} \bm{c}_{1} - \dotsb -
\lambda_{\ell} \bm{c}_{\ell}$ where $\lambda_{1}, \dotsc,
\lambda_{\ell} \in \mathbb{N}$ and $\lambda_{1} + \dotsb +
\lambda_{\ell} = i - 1$.  Therefore, the middle module also vanishes
which proves that $M/g_{j}M$ is weakly $\bm{m}$-regular from level
$k$.
\end{proof}

\begin{theorem} \label{t:regmodule} 
Let $k$ be a positive integer.  If $\scr{C} \subseteq \scr{K}$ and $M$
is weakly $\bm{m}$-regular from level $k$, then $M$ is weakly
$\bm{p}$-regular from level $k$ for every $\bm{p} \in \bm{m} +
\mathbb{N} \scr{C}$.
\end{theorem}

\begin{proof}
Since extension of our base field commutes with the formation of local
cohomology, we may assume for the proof that $\Bbbk$ is infinite.
Because $\scr{C} \subseteq \scr{K}$,
Proposition~\ref{p:nonzerodivisor} implies that we may choose an
almost a nonzerodivisor $g_{j} \in S_{\bm{c}_{j}}$ on $M$ for each $1
\leq j \leq \ell$.

Suppose that $k > 0$ and that $M$ is weakly $\bm{m}$-regular from
level $k$.  Since every $\bm{p} \in \bm{m} + \mathbb{N} \scr{C}$ can
be expressed in the form $\bm{p} = \bm{m} + \lambda_{1} \bm{c}_{1} +
\dotsb + \lambda_{\ell} \bm{c}_{\ell}$ where $\lambda_{1}, \dotsc,
\lambda_{\ell} \in \mathbb{N}$, it suffices to the prove that if $M$
is weakly $\bm{q}$-regular from level $k$ then $M$ is also weakly
$(\bm{q} + \bm{c}_{j})$-regular from level $k$ for each $1 \leq j \leq
\ell$.  We proceed by induction on $\dim M$.  If $\dim M = 0$ then
Grothendieck's vanishing theorem (Theorem~6.1.2 of
\cite{BrodmannSharp}) implies that $H_{B}^{i}(M) = 0$ for $i > 0$.
Thus, $M$ is weakly $\bm{q}$-regular from level $k$ for all $\bm{q}
\in G$ and there is nothing more to prove.  

Assume that $\dim M > 0$ and set $\overline{M} := M/ H_{B}^{0}(M)$.
Since $H_{B}^{0}(M)$ is a $B$-torsion module, the long exact sequence
in cohomology arising from the short exact sequence $0 \longrightarrow
H_{B}^{0}(M) \longrightarrow M \longrightarrow \overline{M}
\longrightarrow 0$ implies that $H_{B}^{i}(M) =
H_{B}^{i}(\overline{M})$ for all $i > 0$.  Hence, $\overline{M}$ is
weakly $\bm{q}$-regular from level $k$ and it suffices to show that
$\overline{M}$ is weakly $(\bm{q} + \bm{c}_{j})$-regular from level
$k$.

Since $g_{j}$ is a nonzerodivisor on $\overline{M}$, we have $\dim
\overline{M}/g_{j} \overline{M} < \dim \overline{M}$.
Lemma~\ref{l:weakly} shows that $\overline{M}/g_{j} \overline{M}$ is
weakly $\bm{q}$-regular from level $k$ and the induction hypothesis
implies $\overline{M}/g_{j} \overline{M}$ is also weakly
$(\bm{q}+\bm{c}_{j})$-regular from level $k$.  Taking the long exact
sequence associated to the exact sequence $0 \longrightarrow
\overline{M}(- \bm{c}_{j}) \longrightarrow \overline{M}
\longrightarrow \overline{M} / g_{j} \overline{M} \longrightarrow 0$,
we obtain the exact sequence $H_{B}^{i}(\overline{M})_{\bm{q}'}
\longrightarrow H_{B}^{i}(\overline{M})_{\bm{q}' + \bm{c}_{j}}
\longrightarrow H_{B}^{i}(\overline{M} / g_{j} \overline{M})_{\bm{q}'
+ \bm{c}_{j}}$.  Since $\overline{M}$ is weakly $\bm{q}$-regular from
level $k$ and $\overline{M}/g_{j} \overline{M}$ is
$(\bm{q}+\bm{c}_{j})$-regular from level $k$, the left and right
modules vanishes when $i \geq k$ and $\bm{q}' = \bm{q} - \lambda_{1}
\bm{c}_{1} - \dotsb - \lambda_{\ell} \bm{c}_{\ell}$ where
$\lambda_{1}, \dotsc, \lambda_{\ell} \in \mathbb{N}$ and $\lambda_{1}
+ \dotsb + \lambda_{\ell} = i - 1$.  Therefore, the middle module also
vanishes which proves $\overline{M}$ is weakly $(\bm{q} +
\bm{c}_{j})$-regular from level $k$.
\end{proof}

Theorem~\ref{t:regmodule} provides the desired alternative
characterization of $\bm{m}$-regularity.

\begin{corollary} \label{c:strongvanishing}
Let $\scr{C} \subseteq \scr{K}$.  The module $M$ is $\bm{m}$-regular
if and only if $M$ is weakly $\bm{m}$-regular from level $1$ and
$H_{B}^{0}(M)_{\bm{p}}=0$ for all $\bm{p} \in \bm{m} + \bigl(
\bigcup_{1 \leq j \leq \ell} (\bm{c}_{j} + \mathbb{N} \scr{C})
\bigr)$.
\end{corollary}

\begin{proof}
Suppose that the module $M$ is $\bm{m}$-regular from level $1$ and
$H_{B}^{0}(M)_{\bm{p}}$ vanishes for all $\bm{p} \in \bm{m} +
\mathbb{N} \scr{C}[1]$.  Since the condition on
$H_{B}^{0}(M)_{\bm{p}}$ is the same as in Definition~\ref{d:regdefn},
we only need to show that $H_{B}^{i}(M)_{\bm{p}} = 0$ for all $\bm{p}
\in \bm{m} + \mathbb{N} \scr{C}[1-i]$.  However, this is the content
of Theorem~\ref{t:regmodule}.  The converse follows directly from the
Definition~\ref{d:regdefn} and Definition~\ref{d:weakreg}.
\end{proof}

The next example illustrates that the condition $\scr{C} \subseteq
\scr{K}$ is necessary for Theorem~\ref{t:regmodule}.  

\begin{example}
Let $S = \Bbbk[x_{1}, x_{2}]$ have the $\mathbb{Z}$-grading defined by
$\deg(x_{1}) = \bm{2}$ and $\deg(x_{2}) = \bm{3}$ and let $B = \langle
x_{1}, x_{2} \rangle$.  Example~\ref{e:standard} shows that $\scr{K} =
6 \mathbb{N}$.  Proposition~\ref{p:firsttopoprop} establishes that the
nonvanishing local cohomology is concentrated in $H_{B}^{2}(S)$ and
that $H_{B}^{2}(S)$ in nonzero exactly in degrees $\bm{-5}, \bm{-7},
\bm{-8}, \bm{-9}, \dotsc$.  If $\scr{C} = \{ \bm{1} \} \not\subset
\scr{K}$, then $S$ is weakly $(\bm{-5})$-regular from level $1$ but
not weakly $(\bm{-4})$-regular from level $1$.  The strategy used in
the proof of Theorem~\ref{t:regmodule} does not apply because
$S_{\bm{1}} = \emptyset$, so there is no nonzerodivisor of degree
$\bm{1}$ .
\end{example}

By design, the regularity of a module measures where its cohomological
complexities vanish.  Since Gr\"{o}bner bases calculations are linked
to homological properties, one expects a strong connection between
regularity and computational complexity.  Both Proposition~\ref{usc}
and Corollary~\ref{c:restoreg} provide further support for this idea
by relating regularity to initial modules and free resolutions.  To
make the connection between regularity and complexity more precise, we
are interested in the following open problem:

\begin{problem}
Give a computationally efficient method of calculating $\reg(M)$.
\end{problem}

When $S$ has the standard grading, this problem is solved in \cite{BS}
by showing that $\reg(M)$ is determined by the largest degree
generator of the initial module $\initial(M)$ with respect to a
reverse lexicographic order in generic coordinates.  Unfortunately,
this technique does not extend directly polynomial rings with
arbitrary multigradings.  As the next example demonstrates, there may
not be any coordinate changes which preserve the grading and change
the initial module.

\begin{example} \label{e:p2blowup}
Let $S$ be the homogeneous coordinate ring of a toric variety $X$
obtained from $\mathbb{P}^{2}$ by a sequence of five blow-ups.  More
explicitly, the minimal lattice points $\scr{B}$ on the rays of the
fan corresponds to the columns of the matrix $\left[
\begin{smallmatrix} 
1 & 1 & 0 & -1 & -1 & -1 & 0 & 1 \\ 0 & 1 & 1 & 1 & 0 & -1 & -1 & -1
\end{smallmatrix} \right]$
and the associated irrelevant ideal is
\begin{align*}
B = \langle&x_3x_4x_5x_6x_7x_8, x_1x_4x_5x_6x_7x_8,
x_1x_2x_5x_6x_7x_8, x_1x_2x_3x_6x_7x_8,\\ 
&x_1x_2x_3x_4x_7x_8, x_1x_2x_3x_4x_5x_8, x_2x_3x_4x_5x_6x_7,
x_1x_2x_3x_4x_5x_6 \rangle \, .
\end{align*}
Hence, $G = \mathbb{Z}^{6}$ and we may assume that $\scr{A}$ be given
by the columns of the matrix
\[
\begin{bmatrix}
1 & 0 & 0 & 0 & 0 & 0 & 1 & -1 \\
0 & 1 & -2 & 1 & 0 & 0 & 0 & 0 \\ 
0 & 0 & 1 & -1 & 1 & 0 & 0 & 0 \\ 
0 & 0 & 0 & 1 & -2 & 1 & 0 & 0 \\ 
0 & 0 & 0 & 0 & 1 & -1 & 1 & 0 \\ 
0 & 0 & 0 & 0 & 0 & 1 & -2 & 1
\end{bmatrix} \, .
\] 
The homogeneous coordinate ring $S = \Bbbk[x_{1}, \dotsc, x_{8}]$ has
the $\mathbb{Z}^{6}$-grading induced by $\deg(x_{i}) = \bm{a}_{i}$ for
$1 \leq i \leq 8$.  Since $S_{\bm{a}_{i}}$ is the $\Bbbk$-span of the
variable $x_{i}$ for $1 \leq i \leq 8$, Corollary~4.7 in \cite{Cox}
establishes that $\Aut(S) = (\Bbbk^{*})^{8}$.  As a consequence, any
change of coordinates which preserves the grading does not alter the
initial ideal.  Hence, one cannot develop a theory of generic initial
ideals.
\end{example}

\section{Degrees of Generators}

The regularity of a module should be regarded as a measure of its
complexity.  In this section, we justify this idea by proving that the
regularity controls the degrees of the minimal generators.  To
understand the minimal generators of a module, we study submodules of
the following form.

\begin{definition}
Let $\scr{D}$ be a subset of $G$.  We define the
\emph{$\scr{D}$-truncation of $M$}, denoted $M|_{\scr{D}}$, to be the
submodule of $M$ generated by all the  homogeneous elements in $M$ of
degree $\bm{p}$ where $\bm{p} \in \scr{D}$.
\end{definition}

In contrast with the standard graded case, the next example
illustrates that $(M|_{\scr{D}})_{\bm{p}}$ may be nonzero even if
$\bm{p} \not\in \scr{D}$.

\begin{example}
Suppose that $S$ is the homogeneous coordinate ring of the Hirzebruch
surface $\mathbb{F}_{2}$; see Examples~\ref{e:hirzebruch} \&
\ref{e:hirzvariety}.  If $\scr{D} := \left[
\begin{smallmatrix} 1 \\ 1 \end{smallmatrix} \right] +
\mathbb{N}^{2}$, then $S|_{\scr{D}}$ is generated in degree $\left[
\begin{smallmatrix} 1 \\ 1 \end{smallmatrix} \right]$.  Observe that
$x_{1}^{3}x_{2}^{2} \in S|_{\scr{D}}$ but $\deg(x_{1}^{3}x_{2}^{2})
\not\in \scr{D}$.
\end{example}

To prove the main theorem in this section, we need the following fact.

\begin{lemma} \label{l:H0}
Let $\bm{c} \in \mathbb{N} \scr{C}$ and let $\scr{V} \subseteq G$ such
that $H_{B}^{0}(M)_{\bm{p}} = 0$ for all $\bm{p} \in \scr{V}$.  If $M$
is $\bm{m}$-regular from level $1$ and $g \in S_{\bm{c}}$ is almost a
nonzerodivisor on $M$, then $H_{B}^{0}(M/g M)_{\bm{p}} = 0$ for all
$\bm{p} \in \scr{V} \cap ( \bm{m} + \bm{c} + \mathbb{N} \scr{C})$.
\end{lemma}

\begin{proof}
Set $\overline{M} := M/ (0 :_{M} g)$.  Because $g$ is almost a
nonzerodivisor on $M$, $(0 :_{M} g)$ is a $B$-torsion module and the
long exact sequence associated to the short exact sequence $0
\longrightarrow (0 :_{M} g) \longrightarrow M \longrightarrow
\overline{M} \longrightarrow 0$ implies that $H_{B}^{i}(M) =
H_{B}^{i}(\overline{M})$ for $i > 0$.  Since $M$ is $\bm{m}$-regular
from level $1$, it follows that $\overline{M}$ is also
$\bm{m}$-regular from level $1$.  Now, the long exact sequence
associated to the exact sequence $0 \longrightarrow \overline{M}(-
\bm{c}) \longrightarrow M \longrightarrow M/gM \longrightarrow 0$
contains $H_{B}^{0}(M)_{\bm{p}} \longrightarrow
H_{B}^{0}(M/gM)_{\bm{p}} \longrightarrow
H_{B}^{1}(\overline{M})_{\bm{p}-\bm{c}}$.  By hypothesis, the left
module vanishes when $\bm{p} \in \scr{V}$.  Since $\overline{M}$ is
$\bm{m}$-regular from level~$1$, the right module vanishes for all
$\bm{p} \in \bm{m} + \bm{c} + \mathbb{N} \scr{C}$.  Hence, the middle
module vanishes for all $\bm{p} \in \scr{V} \cap (\bm{m} + \bm{c} +
\mathbb{N} \scr{C})$.
\end{proof}

We now prove that if $M$ is $\bm{m}$-regular then $M|_{\bm{m}} =
M|_{(\bm{m} + \mathbb{N} \scr{C})}$.  In the standard graded case,
this is true for any $\bm{m} \in \mathbb{Z}$ that is larger than the
maximum degree of the minimal generators.

\begin{theorem} \label{t:gens}
Assume that $\scr{C} \subseteq \scr{K}$.  If the module $M$ is
$\bm{m}$-regular, then we have $M|_{\bm{m}} = M|_{(\bm{m} + \mathbb{N}
\scr{C})}$.  In particular, the degrees of the minimal generators of
$M$ do not belong to the set $\reg(M) + \bigl( \bigcup\nolimits_{1
\leq j \leq \ell} ( \bm{c}_{j} + \mathbb{N} \scr{C}) \bigr)$.
\end{theorem}

\begin{proof}
We prove the following claim: If $M$ is $\bm{m}$-regular from level
$1$ and $H_{B}^{0}(M)_{\bm{p}} = 0$ for all $\bm{p} \in \scr{V}$, then
$M|_{(\bm{m} + \mathbb{N} \scr{C}) \cap \scr{V}}$ is a submodule of
$M|_{\bm{m}}$.  We proceed by induction on $\dim M$.  If $\dim M < 0$,
which is equivalent to saying that $M = 0$, then $M|_{(\bm{m} +
\mathbb{N} \scr{C}) \cap \scr{V}} = 0$ is trivially a submodule of
$M_{\bm{m}}$.

Suppose that $\dim {M} \geq 0$.  Set $\overline{M} := M/ H_{B}^{0}(M)$
and consider the short exact sequence
\begin{equation} \label{gens:1}
0 \longrightarrow H_{B}^{0}(M) \longrightarrow M \longrightarrow
\overline{M} \longrightarrow 0 \, .
\end{equation}  
We claim that it is enough to prove that $\overline{M}|_{(\bm{m} +
\mathbb{N} \scr{C})}$ is a submodule of $\overline{M}|_{\bm{m}}$.  To
see this, let $\bm{p} \in (\bm{m}+\mathbb{N} \scr{C}) \cap \scr{V}$
and let $f \in M_{\bm{p}}$.  If $\overline{f} \in
\overline{M}_{\bm{p}}$ can be written as $\sum_{i} s_{i}
\overline{f}_{i}$ for some $s_{i} \in S$ and some $\overline{f}_{i}
\in \overline{M}_{\bm{m}}$, then $f = \sum_i s_{i} f_{i} + h$ for $h
\in H_{B}^{0}(M)_{\bm{p}}$.  Since $\bm{p} \in \scr{V}$, we have
$H_{B}^{0}(M)_{\bm{p}} = 0$ which implies $h = 0$.  Thus, $f =
\sum_{i} s_{i} f_{i}$ and it suffices to prove the claim for
$\overline{M}$. The long exact sequence associated to \eqref{gens:1}
implies $H_{B}^{i}(M) = H_{B}^{i}(\overline{M})$ for all $i > 0$.
Since $M$ is $\bm{m}$-regular from level $1$ and
$H_{B}^{0}(\overline{M}) = 0$, we deduce that $\overline{M}$ is
$\bm{m}$-regular.  

Extending the base field commutes with the computing local cohomology,
so we may assume without loss of generality that $\Bbbk$ is infinite.
Because $\scr{C} \subseteq \scr{K}$,
Proposition~\ref{p:nonzerodivisor} implies that, for each $1 \leq j
\leq \ell$, we may choose a nonzerodivisor $g_{j} \in S_{\bm{c}_{j}}$
on $\overline{M}$.  It follows from Lemma~\ref{l:weakly} that
$\overline{M}/ g_{j} \overline{M}$ is $\bm{m}$-regular from level $1$
and Lemma~\ref{l:H0} shows that $H_{B}^{0}( \overline{M} / g_{j}
\overline{M})_{\bm{p}} = 0$ for all $\bm{p} \in \bm{m} + \bm{c}_{j} +
\mathbb{N} \scr{C}$.  Since $\dim \overline{M}/ g_{j} \overline{M} <
\dim \overline{M} \leq \dim M$, the induction hypothesis with
$\scr{V}=\bm{m}+\bm{c}_j+\mathbb N \scr{C}$ guarantees that
$(\overline{M}/ g_{j} \overline{M})|_{(\bm{m} + \bm{c}_{j} +
\mathbb{N} \scr{C})}$ is a submodule of $(\overline{M}/ g_{j}
\overline{M})|_{\bm{m}}$.  It remains to show that this implies
$\overline{M}|_{(\bm{m} + \mathbb{N} \scr{C})}$ is a submodule of
$\overline{M}|_{\bm{m}}$.

Suppose otherwise.  For each element $f \in \overline{M}$ with
$\deg(f) \in \bm{m} + \mathbb{N} \scr{C}$, we set
\[
\reach(f) := \min \bigl\{ \lambda_{1} + \dotsb + \lambda_{\ell} :
\text{$\deg(f) = \bm{m} + \lambda_{1} \bm{c}_{1} + \dotsb +
\lambda_{\ell} \bm{c}_{\ell}$ and $\lambda_{1}, \dotsc, \lambda_{\ell}
\in \mathbb{N}$} \bigr\} \, .
\]
Because $\overline{M}$ is noetherian and we are assuming that
$\overline{M}|_{(\bm{m} + \mathbb{N} \scr{C})}
\not\subseteq\overline{M}|_{\bm{m}}$, there is a minimal generator $f
\in \overline{M}|_{(\bm{m} + \mathbb{N} \scr{C})}$ which has smallest
$\reach(f)$ among all the minimal generators of
$\overline{M}|_{(\bm{m} + \mathbb{N} \scr{C})}$ with degree not equal
to $\bm{m}$.  Since $f \not\in \overline{M}|_{\bm{m}}$, we have
$\reach(f) > 0$ and there exists $\bm{c}_{j}$ with $1 \leq j \leq
\ell$ such that $\deg(f) \in \bm{m} + \bm{c}_{j} + \mathbb{N}
\scr{C}$.  From the previous paragraph, we know that the image of $f$
in $(\overline{M}/ g_{j} \overline{M})|_{(\bm{m} + \bm{c}_{j} +
\mathbb{N} \scr{C})}$ belongs to $(\overline{M}/ g_{j}
\overline{M})|_{\bm{m}}$.  Hence, we may choose homogeneous elements
$f_{1}, \dotsc, f_{e} \in \overline{M}|_{\bm{m}}$ and $s_{1}, \dotsc,
s_{e} \in S$ such that $f - s_{1}f_{1} - \dotsb - s_{e}f_{e} \in
(g_{j} \overline{M})|_{(\bm{m} + \bm{c}_{j} + \mathbb{N} \scr{C})}
\subseteq g_{j} \bigl( \overline{M}|_{(\bm{m} + \mathbb{N} \scr{C})}
\bigr)$.  Let $f'$ be the homogeneous element of
$\overline{M}|_{(\bm{m} + \mathbb{N} \scr{C})}$ satisfying the
equation $f - s_{1}f_{1} - \dotsb - s_{e}f_{e} = g_{j} f'$.  Since $f$
does not belong to $\overline{M}|_{\bm{m}}$, the element $f'$ cannot
belong to $\overline{M}|_{\bm{m}}$.  However, $\reach(f') < \reach(f)$
and $f' \in \overline{M}|_{(\bm{m} + \mathbb{N} \scr{C})}$ which
contracts the our choice of $f$.  We conclude that
$\overline{M}|_{(\bm{m} + \mathbb{N} \scr{C})}$ is a submodule of
$\overline{M}_{\bm{m}}$.
\end{proof}

\begin{proof}[Proof of Theorem~\ref{i:gens}]
In the introduction we assumed that $\scr{C}$ was the minimal
generating set of $\scr{K}$.  Hence, Theorem~\ref{i:gens} follows from
Theorem~\ref{t:gens}.
\end{proof}

When $S$ has the standard grading and $\scr{C} = \{ \bm{1} \}$, the
statement that the minimal generators of $M$ do not lie in $\bm{m} +
\mathbb{N} \scr{C}[1]$ is equivalent to saying that the minimal
generators of $M$ have degree at most $\bm{m}$.  In this case,
Theorem~\ref{t:gens} proves that the regularity gives a bound on the
degrees of the minimal generators.  For an ideal in $S$, an upper
bound on the minimal generators yields a finite set containing the
degrees of the minimal generators.  In contrast, Theorem~\ref{t:gens}
does not automatically produce a finite set containing the degrees of
the minimal generators for an ideal when $S$ has a general
multigrading.

We may still ask the question whether if $M$ is $\bm{m}$-regular all
minimal generators of $M$ lie in $\bm{m}-\mathbb N \scr{C}$.  It is
not hard to find examples showing that this need not be the case, but
to date all such examples have had $H^0_B(M) \neq 0$.  In the bigraded
case \cite{HoffmanWang} take a different approach, adding an extra
local cohomology vanishing requirement to guarantee this bound on
generators.

\section{Regularity of $\mathscr{O}_{X}$-modules}

In this section, we develop a multigraded version of regularity for
coherent sheaves on a simplicial toric variety $X$.  Since
$\mathscr{O}_{X}$-modules correspond to finitely generated modules
over the homogeneous coordinate ring of $X$, regularity for
$\mathscr{O}_{X}$-modules is essentially a special case of regularity
for $S$-modules.  Nevertheless, the geometric context provides new
interpretations and applications for regularity.

To begin, we examine the relation between sheaves on a toric variety
$X$ and $G$\nobreakdash-graded modules over its homogeneous coordinate
ring $S$; see \cite{Cox}.  For $\sigma \in \Delta$, we write
$\bm{x}_{\widehat{\sigma}}$ for the monomial $\prod_{i \not\in \sigma}
x_{i} \in S$ and $U_{\sigma} \cong \Spec(S[
\bm{x}_{\widehat{\sigma}}^{-1} ]_{\bm{0}} )$ is the corresponding open
affine subset of $X$.  Every (not necessarily finitely generated)
$G$-graded $S$-module $F$ gives rise to a quasicoherent sheaf
$\widetilde{F}$ on $X$ where $H^{0}(U_{\sigma}, \widetilde{F}) =
(F[\bm{x}_{\widehat{\sigma}}^{-1} ])_{\bm{0}}$.  If $M$ is finitely
generated $G$-graded $S$-module, then $\widetilde{M}$ is a coherent
$\mathscr{O}_{X}$\nobreakdash-module.  Moreover, every quasicoherent
sheaf on $X$ is of this form for some $G$\nobreakdash-graded
$S$-module $F$ and if the sheaf is coherent then $F$ can be taken to
be finitely generated.  For $\bm{p} \in G$, the sheaf associated to
the module $S(\bm{p})$ is denoted $\mathscr{O}_{X}(\bm{p})$.  Note
that the natural map $\mathscr{O}_{X}(\bm{p})
\otimes_{\mathscr{O}_{X}} \mathscr{O}_{X}(\bm{q}) \longrightarrow
\mathscr{O}_{X}(\bm{p} + \bm{q})$ need not be an isomorphism.  For an
$\mathscr{O}_{X}$-module $\mathscr{F}$, we simply write
$\mathscr{F}(\bm{p})$ for $\mathscr{F} \otimes_{\mathscr{O}_{\! X}}
\mspace{-3mu} \mathscr{O}_{X}(\bm{p})$.

The map sending a finitely generated $G$-graded $S$-module $M$ to the
sheaf $\widetilde{M}$ is not injective.  In fact, there are many
nonzero modules which give the zero sheaf.  The following proposition
analyzes this phenomenon.

\begin{proposition} \label{p:zerosheaf}
The sheaf $\widetilde{M}$ is zero if and only if there is $j > 0$ such
that $B^{j} M_{\bm{p}} = 0$ for all $\bm{p} \in \scr{K}$.
\end{proposition}

\begin{proof}
Suppose that $B^{j} M_{\bm{p}} = 0$ for all $\bm{p} \in \scr{K}$.  Fix
$\sigma \in \Delta$ and take $f/\bm{x}^{\bm{u}} \in
(M[\bm{x}_{\widehat{\sigma}}^{-1}])_{\bm{0}}$.
Remark~\ref{r:nonempty} implies that there exists a monomial
$\bm{x}^{\bm{v}} \in S$ such that $\supp(\bm{v}) \subseteq
\widehat{\sigma}$ and $\deg(\bm{x}^{\bm{u} + \bm{v}}) \in
\scr{K}^{\sat}$.  Moreover, we may choose $k \in \mathbb{N}$ such that
$\bm{p} := \deg(\bm{x}^{k \bm{u} + k \bm{v}}) \in \scr{K}$.  It
follows that $\bm{x}^{(k-1)\bm{u} + k \bm{v}} f \in M_{\bm{p}}$ and
$(\bm{x}_{\widehat{\sigma}})^{j} \bm{x}^{(k-1)\bm{u} + k \bm{v}} f =
0$ for some $j > 0$.  Therefore, $f/ \bm{x}^{\bm{u}} = 0$ and
$(M[\bm{x}_{\widehat{\sigma}}^{-1}])_{\bm{0}} = 0$ which shows that
$\widetilde{M}$ is the zero sheaf.

Conversely, suppose that $\widetilde{M}$ is the zero sheaf, fix
$\bm{p} \in \scr{K}$ and let $f \in M_{\bm{p}}$.  Lemma~\ref{l:BandC}
implies that for each $\sigma \in \Delta$ there is a monomial
$\bm{x}^{\bm{u}} \in S$ such that $\deg(\bm{x}^{\bm{u}}) = \bm{p}$ and
$\supp(\bm{u}) \subseteq \widehat{\sigma}$.  The monomial
$\bm{x}^{\bm{u}}$ is invertible in $S[\bm{x}_{\widehat{\sigma}}^{-1}]$
and $f/ \bm{x}^{\bm{u}} \in
(M[\bm{x}_{\widehat{\sigma}}^{-1}])_{\bm{0}} = 0$.  It follows that
$(\bm{x}_{\widehat{\sigma}})^{j} f = 0$ for some $j > 0$.  To see that
one $j$ works for all $f \in M_{\bm{p}}$ and all $\bm{p} \in \scr{K}$,
observe that $M|_{\scr{K}}$ is a submodule of $M$ and thus finitely
generated.  The proposition now follows easily.
\end{proof}

As above, we continue to write $\scr{C} = \{ \bm{c}_{1}, \dotsc,
\bm{c}_{\ell} \}$ for a fixed subset of $G$.  Geometrically, the set
$\scr{C}$ corresponds to choosing a collection of reflexive sheaves on
$X$.  In particular, fixing an ample line bundle on $X$ determines a
set $\scr{C}$ consisting of one element.  At the other extreme, the
toric variety $X$ is equipped with a canonical choice for $\scr{C}$,
namely the minimal generators of the semigroup $\scr{K}$.  We are most
interested in this case.  Nonetheless, we expect that the flexibility
in allowing other sets $\scr{C}$ will also be useful in studying how
maps between toric varieties effect regularity.

In analogy with Definition~\ref{d:regdefn}, we have the following:

\begin{definition} \label{d:sheafreg}
If $\bm{m} \in G$, then an $\mathscr{O}_{X}$-module $\mathscr{F}$ is
\emph{$\bm{m}$-regular} (with respect to $\scr{C}$) if $H^{i} \bigl(
X, \mathscr{F}(\bm{p} ) \bigr) = 0$ for all $i > 0$ and all $\bm{p}
\in \bm{m} + \mathbb{N} \scr{C}[-i]$.  We write $\reg(\mathscr{F})$
for the set $\{ \bm{m} \in G : \text{$\mathscr{F}$ is
$\bm{m}$-regular} \}$.
\end{definition}

\begin{remark} \label{r:mimpliesp}
Notice that if $\mathscr{F}$ is $\bm{m}$-regular then $\mathscr{F}$ is
$\bm{p}$-regular for all $\bm{p} \in \bm{m} + \mathbb{N} \scr{C}$.
\end{remark}

To make an explicit connection between the two definitions of
regularity, we first recall the following.  If $\mathscr{F}$ is the
coherent $\mathscr{O}_{X}$-module corresponding to $M$, then the local
cohomology is related to the (Zariski) cohomology of sheaves by the
exact sequence
\begin{equation} \label{localglobal}
0 \longrightarrow H_{B}^{0}(M) \longrightarrow M \longrightarrow
\bigoplus\limits_{\bm{p} \in G} H^{0} \bigl( X,
\mathscr{F}(\bm{p}) \bigr) \longrightarrow H_{B}^{1}(M)
\longrightarrow 0
\end{equation}
and $H_{B}^{i+1}(M) \cong \bigoplus_{\bm{p} \in G} H^{i} \bigl( X,
\mathscr{F}(\bm{p}) \bigr)$ for all $i \geq 1$; see Proposition~2.3 in
\cite{EMS}.  Using this observation, we have:

\begin{proposition} \label{p:sheafandmodulereg}
Let $\mathscr{F}$ be the coherent $\mathscr{O}_{X}$-module associated
to the module $M$.  The sheaf $\mathscr{F}$ is $\bm{m}$-regular if and
only if $M$ is $\bm{m}$-regular from level~$2$.  The module $M$ is
$\bm{m}$-regular if and only if the following three conditions hold:
\begin{enumerate}
\item $\mathscr{F}$ is $\bm{m}$-regular;
\item the natural map $M_{\bm{p}} \longrightarrow H^{0} \bigl( X,
\mathscr{F}(\bm{p}) \bigr)$ is surjective for all $\bm{p} \in \bm{m} +
\mathbb{N} \scr{C}$; and
\item $H_{B}^{0}(M)_{\bm{p}} = 0$ for all $\bm{p} \in
\bigcup\nolimits_{1 \leq j \leq \ell} \bigl( \bm{m} + \bm{c}_{j} +
\mathbb{N} \scr{C} \bigr)$.
\end{enumerate}
\end{proposition}

\begin{proof}
Since $H_{B}^{i+1}(M)_{\bm{p}} \cong H^{i} \bigl( X,
\mathscr{F}(\bm{p}) \bigr)$ for all $i \geq 1$, we see that
$\mathscr{F}$ is $\bm{m}$-regular if and only if $M$ is
$\bm{m}$-regular from level~$2$.  To establish the second part,
observe that saying $M$ is $\bm{m}$-regular is equivalent to saying
that $M$ is $\bm{m}$-regular from level~2, $H_{B}^{1}(M)_{\bm{p}} = 0$
for all $\bm{p} \in \bm{m} + \mathbb{N} \scr{C}$ and
$H_{B}^{0}(M)_{\bm{p}} = 0$ for all $\bm{p} \in \bm{m} +
\bigcup\nolimits_{1 \leq j \leq \ell} \bigl( \bm{c}_{j} + \mathbb{N}
\scr{C} \bigr)$.  Since the exact sequence \eqref{localglobal} implies
that $H_{B}^{1}(M)_{\bm{p}} = 0$ if and only if $M_{\bm{p}}
\longrightarrow H^{0} \bigl( X, \mathscr{F}(\bm{p}) \bigr)$ is
surjective, the assertion follows.
\end{proof}

\begin{example} \label{e:proj0reg}
Since $H_{B}^{i}(S) = 0$ for $i = 0$, $1$ (see
Remark~\ref{r:whenvanish}), Proposition~\ref{p:sheafandmodulereg}
implies that $\reg(\mathscr{O}_{X}) = \reg(S)$.  In particular,
Example~\ref{e:standardreg} shows that
$\reg(\mathscr{O}_{\mathbb{P}^{d}}) = \mathbb{N}$ when $\scr{C} = \{
\bm{1} \}$.
\end{example}

The next corollary shows that Definition~\ref{d:sheafreg} extends the
original definition given in \S14 in \cite{Mumford}.  When $X =
\mathbb{P}^{d}$ and $\scr{C} = \{ \bm{1} \}$ $\bigl( \bm{1} \in G$
corresponds to $\mathscr{O}_{X}(\bm{1}) \in \Pic(\mathbb{P}^{d})
\bigr)$, the second part of this corollary is Mumford's definition.

\begin{corollary}
Let $\mathscr{F}$ be the coherent $\mathscr{O}_{X}$-module associated
to the module $M$.  If $\scr{C} \subseteq \scr{K}$ then the following
are equivalent:
\begin{enumerate}
\item $\mathscr{F}$ is $\bm{m}$-regular;
\item $H^{i} \bigl( X, \mathscr{F}(\bm{m} - \lambda_{1} \bm{c}_{1} -
\dotsb - \lambda_{\ell} \bm{c}_{\ell} ) \bigr) = 0$ for all $i > 0$
and all $\lambda_{1}, \dotsc, \lambda_{\ell} \in \mathbb{N}$
satisfying $\lambda_{1} + \dotsb + \lambda_{\ell} = i$.
\end{enumerate}
\end{corollary}

\begin{proof}
This follows from Proposition~\ref{p:sheafandmodulereg} and
Corollary~\ref{c:strongvanishing}.
\end{proof}

The regularity of a finite set of points has a geometric
interpretation.  Example~2.16 in \cite{MaclaganSmith2} rephrases the
following result in terms of the associated multigraded Hilbert
polynomial.

\begin{proposition} \label{p:pointsReg}
Assume that $\scr{C} \subseteq \scr{K}$.  Let $Y$ be an artinian
subscheme of $X$ of length $t$ (for example, a set of $t$ points
in $X$) and let $I_{Y}$ be the associated $B$-saturated ideal in $S$.
If $R_{Y} = S/I_{Y}$ then $\bm{m} \in \reg(R_{Y})$ if and only if the
space of forms vanishing on $Y$ has codimension $t$ is the space of
forms of degree $\bm{m}$.
\end{proposition}

\begin{proof}
Since $I_{Y}$ is $B$-saturated, the local cohomology module
$H_{B}^{0}(R_{Y})$ vanishes.  Because $\dim Y = 0$, we also have
$H_{B}^{i+1}(R_{Y}) = \bigoplus_{\bm{p} \in G} H^{i} \bigl( X,
\mathscr{O}_{Y}(\bm{p}) \bigr) = 0$ for all $i > 0$.  Hence,
$\reg(R_{Y})$ is determined by the module $H_{B}^{1}(R_{Y})$.  Since
the exact sequence \eqref{localglobal} becomes $0 \longrightarrow
R_{Y} \longrightarrow \bigoplus_{\bm{p} \in G} H^{0}\bigl( X,
\mathscr{O}_{Y}(\bm{p}) \bigr) \longrightarrow H_{B}^{1}(R_{Y})
\longrightarrow 0$, we see that $\bm{m} \in \reg(R_{Y})$ if and only
if $(R_{Y})_{\bm{m}} \longrightarrow H^{0}\bigl( X,
\mathscr{O}_{Y}(\bm{m}) \bigr)$ is surjective.  Again, because $\dim Y
= 0$, the scheme $Y$ is an affine variety and every reflexive sheaf on
$Y$ is trivial.  It follows that for all $\bm{p}$ we have
$\mathscr{O}_{Y}(\bm{p}) \cong \mathscr{O}_{Y}$ and $\dim_{\Bbbk}
H^{0} \bigl( X, \mathscr{O}_{Y} \bigr) = t$.  Therefore, $R_{Y}$ is
$\bm{m}$-regular if and only if $\dim_{\Bbbk} (R_{Y})_{\bm{m}} = t$.
In other words, the space of forms $(I_{Y})_{\bm{m}}$ that vanish on
$Y$ has codimension $t$.
\end{proof}

Although the definition of regularity may seem rather technical (
``apparently silly'' in Mumford's words), the following theorem
provides a geometric interpretation for regularity.  Concretely, the
regularity of a coherent sheaf measures how much one has to twist for
the sheaf to be generated by its global sections.  We first show that
the module $M$ and certain truncations of $M$ give rise to the same
sheaf.

\begin{lemma} \label{l:trunsheaf}
Assume that $\scr{C} \subseteq \scr{K}$ and $\dim_{\mathbb{R}}
\pos(\overline{\scr{C}}) = r$.  If $\bm{m} \in \mathbb{Z} \scr{K}$ and
$M'$ is the quotient module $M / M|_{(\bm{m} + \mathbb{N} \scr{C})}$
then there exists $j > 0$ such that $B^{j} M_{\bm{p}}' = 0$ for all
$\bm{p} \in \mathbb{Z} \scr{K}$.  Moreover, the module $M$ and its
truncation $M|_{(\bm{m} + \mathbb{N} \scr{C})}$ correspond to the same
$\mathscr{O}_{X}$-module.
\end{lemma} 

\begin{proof} 
Fix $\bm{p} \in \mathbb{Z}\scr{K}$ and $f \in M_{\bm{p}}$.  If $\bm{c}
\in \interior \scr{K}$, then Remark~\ref{r:nonempty} (applied to the
group $\mathbb{Z} \scr{K}$) shows that there is $k \in \mathbb{N}$
such that $\bm{p} + k \bm{c} \in \bm{m} + \scr{K}$.  Since $\scr{C}
\subseteq \scr{K}$ and $\dim_{\mathbb{R}} \pos(\overline{\scr{C}}) =
r$, there exists $\bm{q} \in \scr{K}$ for which $\bm{p} + k \bm{c} +
\bm{q} \in \bm{m} + \mathbb{N}\scr{C}$.  Now, Lemma~\ref{l:BandC}
implies that $B \subseteq \sqrt{ \langle S_{k \bm{c} + \bm{q}}
\rangle}$ which means that there is $j > 0$ such that $B^{j} \subseteq
\langle S_{k \bm{c} + \bm{q}} \rangle$.  Hence, we have $B^{j} f \in
M|_{(\bm{m} + \mathbb{N} \scr{C})}$.  Because $M$ is finitely
generated, there is one $j$ that works for all $f \in M_{\bm{p}}$ and
all $\bm{p} \in \scr{K}$ which establishes the first part.

For the second part, consider the short exact sequence 
\begin{equation} \label{trunses}
0 \longrightarrow M|_{(\bm{m} + \mathbb{N} \scr{C})} \longrightarrow M
\longrightarrow M' \longrightarrow 0 \, .
\end{equation}  
Since the functor $F \mapsto \widetilde{F}$ is exact, it suffices to
show that $\widetilde{M'} = 0$.  By Proposition~\ref{p:zerosheaf},
this is equivalent to the first part of the lemma.
\end{proof}

\begin{theorem} \label{generalMlemma}
Assume that $\scr{C} \subseteq \scr{K}$ and $\dim_{\mathbb{R}}
\pos(\overline{\scr{C}}) = r$.  If $\bm{m} \in \mathbb{Z}\scr{K}$ and
the sheaf $\mathscr{F}$ is $\bm{m}$-regular then the natural map
\begin{equation} \label{generalMlemma:1}
H^{0}\bigl( X, \mathscr{F}(\bm{p}) \bigr) \otimes H^{0}\bigl( X,
\mathscr{O}_{X}(\bm{q}) \bigr) \rightarrow H^{0}\bigl( X,
\mathscr{F}(\bm{p} + \bm{q}) \bigr)
\end{equation}
is surjective for all $\bm{p} \in \bm{m} + \mathbb{N} \scr{C}$ and all
$\bm{q} \in \mathbb{N} \scr{C}$. In particular, the sheaf
$\mathscr{F}(\bm{p})$ is generated by its global sections.
\end{theorem}

\begin{proof}
Consider the $G$-graded $S$-module
\[
M := \Bigl. \Bigl( \bigoplus\nolimits_{\bm{p} \in G} H^{0} \bigl( X,
\mathscr{F}(\bm{p}) \bigr) \Bigr) \Bigr|_{(\bm{m} + \mathbb{N}
\scr{C})} \, .
\]
Because $\mathscr{F}$ is the sheaf associated with the $S$-module
$\bigoplus_{\bm{p} \in \scr{K}} H^{0} \bigl( X, \mathscr{F}(\bm{p})
\bigr)$, Lemma~\ref{l:trunsheaf} guarantees that $\widetilde{M} =
\mathscr{F}$.  Since $\mathscr{F}$ is $\bm{m}$-regular,
Proposition~\ref{p:sheafandmodulereg} shows that $M$ is
$\bm{m}$\nobreakdash-regular from level~$2$.  Combining the definition
of $M$ with the exact sequence \eqref{localglobal}, we have
$H_{B}^{i}(M)_{\bm{p}} = 0$ for all $i = 0$, $1$ and all $\bm{p} \in
\bm{m} + \mathbb{N} \scr{C}$.  Hence, $M$ is $\bm{m}$-regular which
implies that $M$ is also $\bm{p}$-regular for all $\bm{p} \in \bm{m} +
\mathbb{N} \scr{C}$.  If $\bm{p} \in \bm{m} + \mathbb{N} \scr{C}$,
then Theorem~\ref{t:gens} shows that $M|_{\bm{p}} = M|_{(\bm{p} +
\mathbb{N} \scr{C})}$.  In other words, $M_{\bm{p}} \cdot S_{\bm{q}} =
M_{\bm{p}+\bm{q}}$ for all $\bm{q} \in \mathbb{N} \scr{C}$ and this is
equivalent to \eqref{generalMlemma:1} being surjective.  Furthermore,
Lemma~\ref{l:trunsheaf} also guarantees that $\mathscr{F}$ is the
sheaf associated to $M|_{(\bm{p} + \mathbb{N} \scr{C})}$.  Since
$M|_{\bm{p}} = M|_{(\bm{p} + \mathbb{N} \scr{C})}$, this establishes
the last part of the theorem.
\end{proof}

\begin{proof}[Proof of Theorem~\ref{i:sheaves}]
Part~1 follows from Remark~\ref{r:mimpliesp}.  In the introduction,
$\scr{C}$ is the minimal generating set for $\scr{K}$.  Therefore,
Theorem~\ref{generalMlemma} proves Parts~2 and 3.
\end{proof}

When $\mathscr{F} = \mathscr{O}_{X}$, Theorem~\ref{generalMlemma} has
a classical geometric interpretation.  Specifically, if $\bm{p} \in
\reg(\mathscr{O}_{X})$ and $\mathscr{O}_{X}(\bm{p})$ is an ample line
bundle then Theorem~\ref{generalMlemma} implies that the natural map
$H^{0}\bigl( X, \mathscr{O}_{X}(\bm{p}) \bigr) \otimes H^{0}\bigl( X,
\mathscr{O}_{X}(j \bm{p}) \bigr) \longrightarrow H^{0}\bigl( X,
\mathscr{O}_{X}((j+1) \bm{p}) \bigr)$ is surjective for all $j \geq
0$.  It follows (see Exercise~II.5.14 in \cite{hartAG}) that
$\mathscr{O}_{X}(\bm{p})$ embeds $X$ into $\mathbb{P}^{N}$ as a
projectively normal variety where $N := \dim_{\Bbbk} H^{0} \bigl(X,
\mathscr{O}_{X}(\bm{m}) \bigr)$

It is an interesting open problem to describe the set of ample line
bundles on a $d$\nobreakdash-dimensional toric variety $X$ which give
a projectively normal embedding.  When $\mathscr{O}_{X}(\bm{m})$ is an
ample line bundle, \cite{BGT} show that the complete linear series of
the line bundle $\mathscr{O}_{X}\bigl( (d-1) \bm{m} \bigr)$ produces a
projectively normal embedding of $X$.  Moreover, this bound is sharp
on some simplicial toric varieties.  For smooth toric varieties, the
question of whether the map $H^{0} \bigl( X, \mathscr{O}_{X}(\bm{p})
\bigr) \otimes H^{0} \bigl( X, \mathscr{O}_{X}(\bm{q}) \bigr)
\longrightarrow H^{0} \bigl( X, \mathscr{O}_{X}(\bm{p} + \bm{q})
\bigr)$ is surjective for $\mathscr{O}_{X}(\bm{p})$ an ample line
bundle and $\mathscr{O}_{X}(\bm{q})$ a nef line bundle appears in
\cite{Oda}.  This would imply that every ample line bundle induced a
projectively normal embedding.

These questions motivate the study of when $\bm{0} \in
\reg(\mathscr{O}_{X})$.

\begin{proposition} \label{p:productsare0}
If $X$ is a finite product of projective spaces and $\scr{C}$ is the
set of unique minimal generators of $\scr{K}$, then the structure
sheaf $\mathscr{O}_{X}$ is $\bm{0}$-regular.
\end{proposition}

\begin{proof}
We proceed by induction on the number of factors in $X$.  When there
is only one copy of projective space, Example~\ref{e:proj0reg}
establishes the proposition.  For the induction step, we prove the
following stronger claim: Let $X$ be a smooth projective toric variety
such that $\mathscr{O}_{X}$ is $\bm{0}$-regular with respect to $\{
\bm{c}_{1}, \dotsc, \bm{c}_{\ell} \}$.  If $Z = X \times
\mathbb{P}^{h}$ then $\mathscr{O}_{Z}$ is $\bm{0}$-regular with
respect to $\left\{ \left[
\begin{smallmatrix} \bm{c}_{1} \\ \bm{0} \end{smallmatrix}
\right], \dotsc, \left[
\begin{smallmatrix} \bm{c}_{\ell} \\ \bm{0} \end{smallmatrix}
\right], \left[
\begin{smallmatrix} \bm{0} \\ \bm{1} \end{smallmatrix}
\right] \right\}$.

To accomplish this, we first fix some notation.  Let $\Delta$ and
$\Sigma$ be the simplicial complexes arising from the fans of $X$ and
$Z$ respectively.  If the set $\{ \bm{a}_{1}, \dotsc, \bm{a}_{n} \}$
defines the $\mathbb{Z}^{r}$-grading on the homogeneous coordinate
ring of $X$ then the homogeneous coordinate ring $S = \Bbbk[x_{1},
\dotsc, x_{n}, y_{0}, \dotsc, y_{h}]$ has the
$\mathbb{Z}^{r+1}$-grading induced by $\deg(x_{i}) = \left[ 
\begin{smallmatrix} \bm{a}_{i} \\ \bm{0} \end{smallmatrix} 
\right]$ and $\deg(y_{i}) = \left[ 
\begin{smallmatrix} \bm{0} \\ \bm{1} \end{smallmatrix} 
\right]$.  For $\tau \subseteq [n]$ and $\sigma \subseteq \{0, \dotsc,
h\}$, we write $\scr{V}_{\tau \cup \sigma}$ for the subset of
$\mathbb{Z}^{r+1}$ given by
\[
\text{\footnotesize 
$\left[ 
\begin{smallmatrix} 
\sum_{j \in \tau} - \bm{a}_{j} \\ 
- | \sigma | 
\end{smallmatrix} \right]
+ \mathbb{N} \bigl\{  \left[
\begin{smallmatrix} \bm{a}_{j} \\ \bm{0} \end{smallmatrix}
\right] : \text{$j \in \widehat{\tau}$} \bigr\} 
+ \mathbb{N} \bigl\{  \left[
\begin{smallmatrix} -\bm{a}_{j} \\ \bm{0} \end{smallmatrix}
\right] : \text{$j \in \tau$} \bigr\}
+ \mathbb{N} \bigl\{  \left[
\begin{smallmatrix} \bm{0} \\ \bm{1} \end{smallmatrix}
\right] : \text{if $|\sigma| \neq h+1$} \bigr\}
+ \mathbb{N} \bigl\{  \left[
\begin{smallmatrix} \bm{0} \\ \bm{-1} \end{smallmatrix}
\right] : \text{if $|\sigma| \neq 0$} \bigr\}$.}
\]
We also set
\[
\scr{R}^{i} := \left\{ - \lambda_{0} \left[
\begin{smallmatrix} \bm{0} \\ \bm{1} \end{smallmatrix}
\right] -  \lambda_{1} \left[
\begin{smallmatrix} \bm{c}_{1} \\ \bm{0} \end{smallmatrix}
\right] - \dotsb - \lambda_{\ell} \left[
\begin{smallmatrix} \bm{c}_{\ell} \\ \bm{0} \end{smallmatrix}
\right] : \text{$\lambda_{0}, \dotsc, \lambda_{\ell} \in \mathbb{N}$
and $\lambda_{0} + \dotsb + \lambda_{\ell} = i$} \right\} \, .
\]
Now, Corollary~\ref{c:hilbertseries} implies that $\mathscr{O}_{Z}$ is
$\bm{0}$-regular if $\scr{V}_{\tau \cup \sigma} \cap \scr{R}^{i} =
\emptyset$ whenever $i > 0$ and $\widetilde{H}^{i-1}(\Sigma_{\tau \cup
\sigma}) \neq 0$.  We consider four possible cases:
\begin{itemize}
\item If $\emptyset \neq \sigma \neq \{ 0, \dotsc, h \}$, then every
maximal simplex in $\Sigma_{\tau \cup \sigma}$ contains $\sigma$.
Since the induced subcomplex associated to $\sigma$ is the standard
simplex, it follows that $\Sigma_{\tau \cup \sigma}$ is contractible
and $\widetilde{H}^{i-1}(\Sigma_{\tau \cup \sigma}) = 0$.
\item If $\sigma = \emptyset$ then $\Sigma_{\tau \cup \sigma}$ is a
subcomplex of $\Delta$.  Since $\mathscr{O}_{X}$ is $\bm{0}$-regular,
we must have $\scr{V}_{\tau \cup \sigma} \cap \scr{R}^{i} = \emptyset$
whenever $i > 0$ and $\widetilde{H}^{i-1}(\Sigma_{\tau \cup \sigma})
\neq 0$.
\item If $\sigma = \{ 0, \dotsc, h \}$ and $\tau \neq [n]$, then
Alexander duality implies that $\widetilde{H}^{i-1}(\Sigma_{\tau \cup
\sigma})$ is isomorphic to $\widetilde{H}^{d+h-i-1}
(\Sigma_{\widehat{\tau \cup \sigma}})$ which reduces us to the second
case.
\item If $\sigma = \{ 0, \dotsc, h \}$ and $\tau = [n]$ then
$\widetilde{H}^{i-1}(\Sigma_{\tau \cup \sigma}) \neq 0$ if and only if
$i = d + h$.  Since both $\mathscr{O}_{\mathbb{P}^{d}}$ and
$\mathscr{O}_{X}$ are $\bm{0}$-regular, we must have $\scr{V}_{\tau
\cup \sigma} \cap \scr{R}^{d+d'} = \emptyset$. \qedhere
\end{itemize}
\end{proof}

The following example demonstrates that even when $X$ is a smooth Fano
toric variety the structure sheaf $\mathscr{O}_{X}$ may not be
$\bm{0}$-regular.

\begin{example}
Let $X$ be a blowup of $\mathbb{P}^{2}$ at three points.
Specifically, the lattice points $\scr{B}$ correspond to the columns
of the matrix $\left[
\begin{smallmatrix}
1 & 0 & -1 & -1 & 0 & 1 \\
0 & 1 & 1 & 0 & -1 & -1
\end{smallmatrix} 
\right]$.  The regular triangulation of $\scr{B}$ is encoded in the
irrelevant ideal 
\[
B = \langle x_{1}x_{2}x_{3}x_{4},
x_{2}x_{3}x_{4}x_{5}, x_{3}x_{4}x_{5}x_{6}, x_{1}x_{4}x_{5}x_{6},
x_{1}x_{2}x_{5}x_{6}, x_{1}x_{2}x_{3}x_{6} \rangle \, .
\]  
We have $G = \mathbb{Z}^{4}$ and we may assume that $\scr{A}$ is given
by the columns of the matrix
\[
\begin{bmatrix}
1 & 0 & 0 & 0 & 1 & -1 \\
0 & 1 & -1 & 1 & 0 & 0 \\
0 & 0 & 1 & -1 & 1 & 0 \\
0 & 0 & 0 & 1 & -1 & 1
\end{bmatrix} \, .
\]
Let $\scr{C}$ be the unique minimal generators of the monoid $\scr{K}
= \scr{K}^{\sat}$:
\[
\scr{C} = \left\{ \left[
\begin{smallmatrix} 1 \\ 1 \\ 0 \\ 0 \end{smallmatrix}
\right], \left[
\begin{smallmatrix} 0 \\ 0 \\ 1 \\ 0 \end{smallmatrix}
\right], \left[
\begin{smallmatrix} 1 \\ 0 \\ 1 \\ 0 \end{smallmatrix}
\right], \left[
\begin{smallmatrix} 0 \\ 0 \\ 0 \\ 1 \end{smallmatrix}
\right], \left[
\begin{smallmatrix} 0 \\ 1 \\ 0 \\ 1 \end{smallmatrix}
\right] \right\} \, .
\]
It follows from Corollary~\ref{c:hilbertseries} that
\[
H^{2} \bigl( X, \mathscr{O}_{X}(- \bm{c}_{3} - \bm{c}_{5}) \bigr) =
H^{2} \bigl( X, \mathscr{O}_{X}(- \bm{a}_{1} - \bm{a}_{2} - \bm{a}_{3}
- \bm{a}_{4} - \bm{a}_{5} - \bm{a}_{6}) \bigr) =
\widetilde{H}^{1}(\Delta) \neq 0
\]
which implies $\bm{0} \not\in \reg(\mathscr{O}_{X})$.
\end{example}

The following remains an interesting open problem.

\begin{problem}
Let $\scr{C}$ be a set of minimal generators of $\scr{K}$.  Give a
combinatorial characterization of all toric varieties $X$ such that
$\scr{K} \subseteq \reg(\mathscr{O}_{X})$.
\end{problem}

We end this section with an upper semi-continuity result.  This is
well-known when $X = \mathbb{P}^{d}$.  We write $\initial(I)$ for the
initial ideal of an ideal $I \subseteq S$ with respect to some
monomial order.

\begin{proposition} \label{usc}
Let $S$ be the homogeneous coordinate ring of a simplicial toric
variety.  If $I$ is an ideal in $S$, then $\reg \bigl( S/\initial(I)
\bigr) \subseteq \reg(S/I)$.  Moreover, if $I$ is
$B$\nobreakdash-saturated and $J = \bigl( \initial(I) : B^{\infty}
\bigr)$ then $\reg(S/J) \subseteq \reg(S/I)$.
\end{proposition}

\begin{proof}
Fix $j > 0$. If $\bm{x}^{\bm{u}} f$ belongs to $I$ for all
$\bm{x}^{\bm{u}} \in B^{j}$, then $\bm{x}^{\bm{u}} \initial(f)$
belongs to $\initial(I)$.  Hence, $f \in (I:B^{\infty})$ implies that
$\initial(f) \in ( \initial(I) : B^{\infty})$.  We conclude that
$\dim_{\Bbbk} H_{B}^{0}(S/I)_{\bm{p}} \leq \dim_{\Bbbk} H_{B}^{0}( S/
\initial(I))_{\bm{p}}$ for all $\bm{p} \in G$.  Theorem~15.17 in
\cite{E} gives a flat family over $\mathbb{A}^{1}$ whose general fiber
is $S/I$ and whose special fiber is $S/(\initial(I))$.  Since
$H_{B}^{i+1}(F) \cong \bigoplus_{\bm{p} \in G} H^{i}\big(X,
\mathscr{F}(\bm{p}) \big)$ for all $i \geq 1$, Theorem~12.8 in
\cite{hartAG} yields $\dim_{\Bbbk} H_{B}^{i} ( S/ I )_{\bm{p}} \leq
\dim_{\Bbbk} H_{B}^{i} ( S/ \initial(I) )_{\bm{p}}$ for all $i > 1$.
If $\mathscr{I}$ denotes the sheaf of ideals corresponding to $I$,
then the exact sequence \eqref{localglobal} gives
\[
\dim_{\Bbbk} H_{B}^{0}(S/I)_{\bm{p}} + \dim_{\Bbbk} H^{0} \bigl(
X,(\mathscr{O}_{X}/\mathscr{I})(\bm{p}) \bigr) = \dim_{\Bbbk}
(S/I)_{\bm{p}} +\dim_{\Bbbk} H_{B}^{1}(S/I)_{\bm{p}} \, .
\]
Since the terms on left-hand side do not decrease and $\dim_{\Bbbk}
(S/I)_{\bm{p}}$ is constant when passing to an initial ideal, we see
that $\dim_{\Bbbk} H_{B}^{1}(S/I)_{\bm{p}}$ also does not decrease
when passing to an initial ideal.  Thus, we conclude that $H_{B}^{i} (
S/ \initial(I))_{\bm{p}} = 0$ implies $H_{B}^{i}( S/I )_{\bm{p}} = 0$
and the first assertion follows.  Since $H_{B}^{0}(S/J) = 0$ and
$H_{B}^{i}(S/J) = H_{B}^{i} \bigr( S/\initial(I) \bigr)$ for $i > 0$,
this also establishes the second assertion.
\end{proof}

\section{Syzygies and Chain Complexes} \label{complexsection}

The final section of the paper examines the relationship between the
regularity of $M$ and the syzygies of $M$.  We give a combinatorial
formula, involving only the degrees of the syzygies and the regularity
of $S$, for a subset of $\reg(M)$.  When $S$ is the homogeneous
coordinate ring of weighted projective space, this subset actually
equals $\reg(M)$.  In other words, we recover the characterization of
regularity in terms of Betti numbers in this case; see Theorem~5.5 in
\cite{Benson}.  We extend this formula  to coherent
$\mathscr{O}_{X}$-modules.  Changing directions, we then describe
chain complexes associated to certain elements of $\reg(M)$.  At the
level of $\mathscr{O}_{X}$-modules, we obtain a locally free
resolution of $\mathscr{F}$ from specific elements in
$\reg(\mathscr{F})$.

If $\mathbf{E}$ is a free resolution of the module $M$, then there is
a spectral sequence that computes the local cohomology of $M$ from
$\mathbf{E}$, namely $E_{2}^{i,j} := H^{i} \bigl(
H_{B}^{j}(\mathbf{E}) \bigr) \Longrightarrow H_{B}^{i+j}(M)$.  When
$S$ is the homogeneous coordinate ring of weighted projective space,
this spectral sequence degenerates because $H_{B}^{i}(S)$ vanishes
unless $i = d+1$; see Remark~\ref{r:whenvanish}.  As a consequence,
there is a simple characterization of the regularity of $M$ in terms
of the degrees of the syzygies.  In contrast,
Remark~\ref{r:whenvanish} also shows that nonvanishing local
cohomology of $S$ is not typically concentrated in a single
cohomological degree.  Hence, the spectral sequence does not
degenerate and one cannot expect as simple a relationship between
regularity and syzygies in the general situation.  Despite this, the
syzygies of $M$ do provide a method for approximating the regularity
of $M$ which captures the description in the special case.  Moreover,
this technique also works on a larger class of chain complexes
associated to $M$.

We start by describing how regularity behaves in short exact
sequences.  

\begin{lemma} \label{l:ses}
If $0 \longrightarrow M' \longrightarrow M \longrightarrow M''
\longrightarrow 0$ is a short exact sequence of finitely generated
$G$-graded $S$-modules, then we have the following:
\begin{enumerate}
\item $\reg^{i}(M') \cap \reg^{i}(M'') \subseteq \reg^{i}(M)$,
\item $\bigl( \bigcup\nolimits_{1 \leq j \leq \ell} \bigl( -
\bm{c}_{j} +\reg^{i+1}(M')\bigr) \bigr) \cap \reg^{i}(M) \subseteq
\reg^{i}(M'')$ and
\item $\reg^{i}(M) \cap \bigl( \bigcap\nolimits_{1 \leq j \leq \ell}
\bigl( \bm{c}_{j} + \reg^{i-1}(M'') \bigr) \bigr) \subseteq
\reg^{i}(M')$.
\end{enumerate}
\end{lemma}

\begin{proof}
The associated long exact sequence in cohomology contains the exact
sequence
\begin{equation} \label{ses:1}
H_{B}^{i-1}(M'')_{\bm{p}} \longrightarrow H_{B}^{i}(M')_{\bm{p}}
\longrightarrow H_{B}^{i}(M)_{\bm{p}} \longrightarrow
H_{B}^{i}(M'')_{\bm{p}} \longrightarrow H_{B}^{i+1}(M')_{\bm{p}} \, .
\end{equation}
Suppose that both $M'$ and $M''$ are $\bm{m}$-regular from level $k$.
This means that the second and fourth modules in \eqref{ses:1} vanish
for all $i \geq k$ and all $\bm{p} \in \bm{m} + \mathbb{N}
\scr{C}[1-i]$.  Hence, the third module in \eqref{ses:1} also vanishes
for all $i \geq k$ and all $\bm{p} \in \bm{m} + \mathbb{N}
\scr{C}[1-i]$ which implies that $M$ is $\bm{m}$-regular from level
$k$.

Similarly, if $M$ is $\bm{m}$-regular from level $k$ and $M'$ is
$(\bm{m} + \bm{c}_{j})$-regular from level $(k+1)$ for some $1 \leq j
\leq \ell$, then the third and fifth modules in \eqref{ses:1} vanish
for all $i \geq k$ and all $\bm{p} \in \bm{m} + \mathbb{N}
\scr{C}[1-i]$.  It follows that the fourth module in \eqref{ses:1}
also vanishes under the same conditions.  Thus, $M''$ is
$\bm{m}$-regular from level $k$.

Finally, if $M$ is $\bm{m}$-regular from level $k$ and $M''$ is
$(\bm{m} - \bm{c}_{j})$-regular from level $(k-1)$ for every $1 \leq j
\leq \ell$, then the first and third modules in \eqref{ses:1} vanish
for all $i \geq 0$ and all $\bm{p} \in \bm{m} + \mathbb{N}
\scr{C}[1-i]$.  Hence, the second module in \eqref{ses:1} vanishes
under the same conditions and therefore $\bm{m} \in \reg^{k}(M')$.
\end{proof}

The difference between the union in Part~2 and the intersection in
Part~3 introduces an asymmetry in working with short exact sequences.
In many cases, this prevents one from giving a simple characterization
of $\bm{m}$-regularity in terms of minimal free resolutions.  However,
when $\scr{C}$ consists of a single element, as in the standard graded
case, this obstruction is not present.  In this case, we have the more
symmetric: $\bigl(- \bm{c}_{1} +\reg^{i+1}(M') \bigr) \cap \reg^{i}(M)
\subseteq \reg^{i}(M)$ and $\bigl(\bm{c}_{1} +\reg^{i-1}(M'') \bigr)
\cap \reg^{i}(M) \subseteq \reg^{i}(M')$.

The approach used in the proof of Lemma~\ref{l:ses} also leads to an
analogous result for a short exact sequence of coherent
$\mathscr{O}_{X}$-modules.  However, there is a notable change in
Part~3.  Specifically, we must assume that $i > 1$ because the
hypothesis that $\mathscr{F}$ is $\bm{m}$-regular does not place any
conditions on $H^{0} \bigl( X, \mathscr{F}(\bm{p}) \bigr)$.

The next theorem provides a method for estimating the regularity of a
module $M$ from certain chain complexes associated to $M$.  We say a
chain complex of $S$-modules has \emph{$B$-torsion homology} if every
homology module is a $B$-torsion module.

\begin{theorem} \label{t:seqtoreg}
Let $\mathbf{E} := \dotsb \longrightarrow E_{3} \xrightarrow{\;\;
\partial_{3} \;\;} E_{2} \xrightarrow{\;\; \partial_{2} \;\;} E_{1}
\xrightarrow{\;\; \partial_{1} \;\;} E_{0} \longrightarrow 0$ be a
chain complex of finitely generated $G$-graded $S$-modules with
$B$-torsion homology.  If $\partial_{0} \colon E_{0} \longrightarrow
M$ is a surjective map then we have
\[
\bigcup_{\phi \colon [d+1] \rightarrow [\ell]} \left( \bigcap_{0 \leq
i \leq d+1} \bigl( - \bm{c}_{\phi(1)} - \dotsb - \bm{c}_{\phi(i)} +
\reg^{i}(E_{i}) \bigr) \right) \subseteq \reg(M), \, 
\]
where the union is over all functions $\phi \colon [d+1]
\longrightarrow [\ell]$.
\end{theorem}

\begin{proof}
Fix a function $\phi \colon [d+1] \longrightarrow [\ell]$.  We claim
that for all $k \geq 0$
\[
\bigcap_{k \leq i \leq d+1} \bigl( - \bm{c}_{\phi(k+1)} - \dotsb -
\bm{c}_{\phi(i)} + \reg^{i}(E_{i}) \bigr) \subseteq \reg^{k}(\image
\partial_{k}) \, .
\]  
Since $\image \partial_{0} = M$, this will prove the theorem.  We
establish the claim by using a descending induction on $k$.  Since
$H_{B}^{j}(\image \partial_{k})$ vanishes for all $j > d+1$, we have
$\reg^{k}(\image \partial_{k}) = G$ for $k > d+1$ and the claim holds.
Suppose $k \leq d+1$ and consider the exact sequence $0
\longrightarrow H_{k} (\mathbf{E}) \longrightarrow E_{k}/ \image
\partial_{k+1} \longrightarrow \image \partial_{k} \longrightarrow 0$.
Because $H_{k}(\mathbf{E})$ is a $B$-torsion module, the associated
long exact sequence in cohomology implies that $H_{B}^{j}(E_{k}/
\image \partial_{k+1}) = H_{B}^{j}(\image \partial_{k})$ for all $j >
0$ and that $H_{B}^{0}(E_{k}/ \image \partial_{k+1})$ surjects onto
$H_{B}^{0}(\image \partial_{k})$.  It follows that $\reg^{k}(E_{k}/
\image \partial_{k+1}) \subseteq \reg^{k}(\image \partial_{k})$ for
all $k \geq 0$.  Applying Lemma~\ref{l:ses} to the short exact
sequence $0 \longrightarrow \image \partial_{k+1} \longrightarrow
E_{k} \longrightarrow E_{k} / \image \partial_{k+1} \longrightarrow 0$
gives $\reg^{k+1}(\image \partial_{k+1})[-1] \cap \reg^{k}(E_{k})
\subseteq \reg^{k}(E_{k} / \image \partial_{k+1})$.  Since the
induction hypothesis implies that $\bigcap_{k+1 \leq i \leq d+1}
\bigl( - \bm{c}_{\phi(k+2)} - \dotsb - \bm{c}_{\phi(i)} +
\reg^{i}(E_{i}) \bigr) \subseteq \reg^{k+1}(\image \partial_{k+1})$,
we have $\bigcap_{k+1 \leq i \leq d+1} \bigl( - \bm{c}_{\phi(k+1)} -
\dotsb - \bm{c}_{\phi(i)} + \reg^{i}(E_{i}) \bigr) \subseteq
\reg^{k+1}(\image \partial_{k+1})[-1]$.  Combining these three
inclusions, we obtain $\bigcap_{k \leq i \leq d+1} \bigl( -
\bm{c}_{\phi(k+1)} - \dotsb - \bm{c}_{\phi(i)} + \reg^{i}(E_{i})
\bigr) \subseteq \reg^{k}(\image \partial_{k})$.
\end{proof}

Since the proof of Theorem~\ref{t:seqtoreg} only used Part~2 of
Lemma~\ref{l:ses}, the same approach leads to a formula for estimating
the regularity of a sheaf $\mathscr{F}$ from a resolution.  We leave
the details to the interested reader.  Because it is frequently used
in applications, we explicitly state the following special case of
Theorem~\ref{t:seqtoreg}.

\begin{corollary} \label{c:restoreg}
If $0 \longrightarrow E_{s} \longrightarrow \dotsb \longrightarrow
E_{2} \longrightarrow E_{1} \longrightarrow E_{0}$ is a free
resolution of the module $M$ with $E_{i} = \bigoplus\limits_{1 \leq j
\leq h_{i}} S( - \bm{q}_{i,j})$ for some $\bm{q}_{i,j} \in G$, then we
have
\begin{equation} \label{restoreg:1}
\bigcup_{\phi \colon [d+1] \rightarrow [\ell]} \Biggl(
\bigcap_{\begin{subarray}{c} 
0 \leq i \leq \min\{ d+1, s \} \\ 
1 \leq j \leq h_{i} 
\end{subarray}} 
\bigl( \bm{q}_{i,j} - \bm{c}_{\phi(1)} - \dotsb - \bm{c}_{\phi(i)} +
\reg^{i}(S) \bigr) \Biggr) \subseteq \reg(M) \, .
\end{equation}
\end{corollary}

\begin{proof}
Since $\reg^{i}(E_{i}) = \bigcap_{1 \leq j \leq h_{i}} \bigl(
\bm{q}_{i,j} + \reg^{i}(S) \bigr)$, \eqref{restoreg:1} follows from
Theorem~\ref{t:seqtoreg}
\end{proof}

The next three examples illustrate this corollary.

\begin{example} \label{e:standardsyzygies}
Suppose $S$ is the homogeneous coordinate ring of weighted projective
space; see Example~\ref{e:standard}.  Let $\scr{C} = \{ \bm{1} \}$ and
let $0 \longrightarrow E_{s} \longrightarrow \dotsb \longrightarrow
E_{2} \longrightarrow E_{1} \longrightarrow E_{0}$ be the minimal free
resolution of the module $M$ where $E_{i} = \bigoplus\nolimits_{1 \leq
j \leq h_{i}} S( - \bm{q}_{i,j})$ for some $\bm{q}_{i,j} \in
\mathbb{Z}$.  Hilbert's syzygy theorem implies that $s \leq n = d+1$.
Hence, Example~\ref{e:standardreg} and \eqref{restoreg:1} imply that
$\max_{i,j} \{ \bm{q}_{i,j} - i + n - \bm{a}_{1} - \dotsb - \bm{a}_{n}
\} \in \reg(M)$.  In particular, if $S$ has the standard grading, then
$M$ is $(\max_{i,j} \{ \bm{q}_{i,j} - i \})$-regular.
\end{example}

\begin{example} \label{e:basicRes}
Let $S$ be the homogeneous coordinate ring of 
$\mathbb{P}^{1} \times \mathbb{P}^{1} \times \mathbb{P}^{1}$.  This
means that $\scr{B}$ corresponds to the columns of the matrix
\[
\begin{bmatrix}
1 & -1 & 0 & 0 & 0 & 0 \\ 
0 & 0 & 1 & -1 & 0 & 0 \\
0 & 0 & 0 & 0 & 1 & -1
\end{bmatrix} \, .
\]
and the irrelevant ideal is $B = \langle x_{1}, x_{2} \rangle \cap
\langle x_{3}, x_{4} \rangle \cap \langle x_{5}, x_{6} \rangle$.
Hence, $G = \mathbb{Z}^{3}$ and we may assume that $\scr{A}$ is given
by the columns of the matrix
\[
\begin{bmatrix} 
1 & 1 & 0 & 0 & 0 & 0 \\ 
0 & 0 & 1 & 1 & 0 & 0 \\
0 & 0 & 0 & 0 & 1 & 1
\end{bmatrix} \, .
\]
The polynomial ring $S = \Bbbk[x_{1}, \dotsc, x_{6}]$ has
the corresponding $\mathbb{Z}^{3}$-grading and $\scr{K} =
\scr{K}^{\sat} = \mathbb{N}^{3}$.  If $\scr{C} = \{ \bm{c}_{1},
\bm{c}_{2}, \bm{c}_{3} \} \subset \scr{K}$ where $\bm{c}_{i}$ is the
$i$th standard basis vector in $\mathbb{Z}^{3}$, then
Corollary~\ref{c:hilbertseries}) implies that $\reg(S) =
\mathbb{N}^{3}$.

Consider the module $M = S/I$ where $I = \langle x_{1} - x_{2}, x_{3}
- x_{4} , x_{5} - x_{6} \rangle$.  Since $I$ is a reduced prime ideal
of codimension $3$, the subvariety of $\mathbb{P}^{1} \times
\mathbb{P}^{1} \times \mathbb{P}^{1}$ defined by $I$ is a single
point.  The minimal free resolution of $M$ has the form
\[
0 \longrightarrow S(- \bm{c}_{1} - \bm{c}_{2} - \bm{c}_{3})
\longrightarrow \bigoplus_{1 \leq i < j \leq 3} S(- \bm{c}_{i} -
\bm{c}_{j}) \longrightarrow \bigoplus_{1 \leq i \leq 3} S(-\bm{c}_{i})
\longrightarrow S \, .
\]
A calculation using Corollary~\ref{c:restoreg} shows that $\bigcup_{1
\leq i < j \leq 3} ( \bm{c}_{i} + \bm{c}_{j} + \mathbb{N}^{3})
\subseteq \reg(M)$.  Note that this is not the entire regularity of
$M$.  Proposition~\ref{p:pointsReg} show that $\reg(M) =
\mathbb{N}^{3}$.
\end{example}

\begin{example}
As in Example~\ref{e:basicRes}, we assume that $S$ is the homogeneous
coordinate ring of the toric variety $\mathbb{P}^{1} \times
\mathbb{P}^{1} \times \mathbb{P}^{1}$ and we choose $\scr{C}$ to be
the set of standard basis vectors in $\mathbb{Z}^{3}$.  Consider $M =
S/ I$ where $I = \langle x_{1} - x_{2}, x_{3} - x_{4} , x_{5} - x_{6}
\rangle \cap \langle x_{1} - 2x_{2}, x_{3} - 2x_{4} , x_{5} - 2x_{6}
\rangle$.  The $B$-saturated ideal $I$ corresponds to two distinct
points on $\mathbb{P}^{1} \times \mathbb{P}^{1} \times
\mathbb{P}^{1}$.  The minimal free resolution of $M$ has the form
\[
0 \rightarrow \bigoplus_{\bm{q} \in \scr{D}_{5}} S(- \bm{q})
\rightarrow \bigoplus_{\bm{q} \in \scr{D}_{4}} S(- \bm{q})
\rightarrow \bigoplus_{\bm{q} \in \scr{D}_{3}} S(- \bm{q})
\rightarrow \bigoplus_{\bm{q} \in \scr{D}_{2}} S(- \bm{q})
\rightarrow \bigoplus_{\bm{q} \in \scr{D}_{1}} S(- \bm{q})
\rightarrow S
\]
where, for $1 \leq i \leq 5$, the set $\scr{D}_{i}$ is given by the
column vectors of the matrix $D_{i}$ and
\begin{align*}
D_{1} &= \left[ \begin{smallmatrix}
2 & 1 & 1 & 1 & 1 & 0 & 0 & 0 & 0 \\
0 & 1 & 1 & 0 & 0 & 2 & 1 & 1 & 0 \\
0 & 0 & 0 & 1 & 1 & 0 & 1 & 1 & 2
\end{smallmatrix} \right] \, ,
&D_{2} &= \left[ \begin{smallmatrix}
2 & 2 & 2 & 2 & 1 & 1 & 1 & 1 & 1 & 1 & 1 & 1 & 1 & 1 & 0 & 0 & 0 & 0
\\
1 & 1 & 0 & 0 & 2 & 2 & 1 & 1 & 1 & 1 & 1 & 1 & 0 & 0 & 2 & 2 & 1 & 1
\\
0 & 0 & 1 & 1 & 0 & 0 & 1 & 1 & 1 & 1 & 1 & 1 & 2 & 2 & 1 & 1 & 2 & 2 
\end{smallmatrix} \right] \, , \\
D_{3} &= \left[ \begin{smallmatrix}
2 & 2 & 2 & 2 & 2 & 2 & 1 & 1 & 1 & 1 & 1 & 1 & 1 & 1 & 0 \\
2 & 1 & 1 & 1 & 1 & 0 & 2 & 2 & 2 & 2 & 1 & 1 & 1 & 1 & 2 \\
0 & 1 & 1 & 1 & 1 & 2 & 1 & 1 & 1 & 1 & 2 & 2 & 2 & 2 & 2
\end{smallmatrix} \right] \, , 
&D_{4} &= \left[ \begin{smallmatrix}
2 & 2 & 2 & 2 & 1 & 1 \\ 
2 & 2 & 1 & 1 & 2 & 2 \\ 
1 & 1 & 2 & 2 & 2 & 2 
\end{smallmatrix} \right] \, , \qquad
D_{5} = \left[ \begin{smallmatrix}
2 \\ 2 \\ 2
\end{smallmatrix} \right] \, .
\end{align*}
Observe that the length of the minimal free resolution is greater than
$d+1$.  Applying Corollary~\ref{c:restoreg}, we deduce that
\[
\left( \left[
\begin{smallmatrix} 1 \\ 2 \\ 2 \end{smallmatrix} \right] +
\mathbb{N}^{3} \right) \cup \left( \left[ 
\begin{smallmatrix} 2 \\ 1 \\ 2 \end{smallmatrix} \right] +
\mathbb{N}^{3} \right) \cup \left( \left[ 
\begin{smallmatrix} 2 \\ 2 \\ 1 \end{smallmatrix} \right] +
\mathbb{N}^{3} \right) \subset \reg(M) \, .
\]
However, Proposition~\ref{p:pointsReg} shows that 
\[
\reg(M) = \left( \left[
\begin{smallmatrix} 1 \\ 0 \\ 0 \end{smallmatrix} \right] +
\mathbb{N}^{3} \right) \cup \left( \left[ 
\begin{smallmatrix} 0 \\ 1 \\ 0 \end{smallmatrix} \right] +
\mathbb{N}^{3} \right) \cup \left( \left[ 
\begin{smallmatrix} 0 \\ 0 \\ 1 \end{smallmatrix} \right] +
\mathbb{N}^{3} \right) \, .
\]
\end{example}

To describe a partial converse for Theorem~\ref{t:seqtoreg}, we
concentrate on the situation arising from smooth toric varieties.
Specifically, we assume that $\scr{K} = \scr{K}^{\sat}$ and that
$\scr{C}$ is a finite generating set for the monoid $\scr{K}$.

\begin{lemma} \label{l:trunreg}
Let $\mathbb{N} \scr{C} = \scr{K} = \scr{K}^{\sat}$.  If $M$ is
$\bm{m}$-regular then $M|_{(\bm{m} + \scr{K})}$ is also
$\bm{m}$-regular.
\end{lemma}

\begin{proof}
Suppose that $M$ is $\bm{m}$-regular and let $M'$ be the quotient
module $M / M|_{(\bm{m} + \scr{K})}$.  Since $\mathbb{N} \scr{C} =
\scr{K} = \scr{K}^{\sat}$, Lemma~\ref{l:trunsheaf} establishes that
$M'$ is a $B$-torsion module.  Hence, the long exact sequence
associated to \eqref{trunses} implies that $H_{B}^{i}(M|_{(\bm{m} +
\scr{K})}) = H_{B}^{i}(M)$ for all $i > 1$ and $H_{B}^{1}(M|_{(\bm{m}
+ \scr{K})})_{\bm{p}} = H_{B}^{1}(M)_{\bm{p}}$ for all $\bm{p} \in
\bm{m} + \scr{K}$.  For $k > 0$, we deduce that $M$ is
$\bm{m}$-regular from level $k$ if and only if the submodule
$M|_{(\bm{m} + \scr{K})}$ is $\bm{m}$-regular from level $k$.  The
long exact sequence also shows that $H_{B}^{0}(M|_{(\bm{m} +
\scr{K})})$ is a submodule of $H_{B}^{0}(M)$.  Thus, if
$H_{B}^{0}(M)_{\bm{p}} = 0$ then $H_{B}^{0}(M|_{(\bm{m} +
\scr{K})})_{\bm{p}}$ also vanishes.  Therefore, $M|_{(\bm{m} +
\scr{K})}$ is $\bm{m}$-regular.
\end{proof}

\begin{theorem} \label{t:existscomplex}
Assume that $\mathbb{N} \scr{C} = \scr{K} = \scr{K}^{\sat}$.  If
$\bm{c} \in \reg(S) \cap \bigl( \bigcap_{1 \leq j \leq \ell}
(\bm{c}_{j} + \scr{K}) \bigr)$ and $\bm{m} \in \reg(M)$ then there
exists
\begin{enumerate}
\item a chain complex $\;\; \dotsb \longrightarrow E_{3}
\xrightarrow{\;\; \partial_{3} \;\;} E_{2} \xrightarrow{\;\;
\partial_{2} \;\;} E_{1} \xrightarrow{\;\; \partial_{1} \;\;} E_{0}$
with $B$-torsion homology and $E_{i} = \bigoplus_{1 \leq j \leq h_{i}}
S(- \bm{m} - i \bm{c})$, and
\item a surjective map $\partial_{0} \colon E_{0} \longrightarrow
M|_{(\bm{m} + \mathbb{N} \scr{C})}$.
\end{enumerate}
\end{theorem}

\begin{proof}
We prove by induction on $k$ that there exists a chain complex with
$B$-torsion homology
\begin{equation} \label{existscomplex:1}
0 \longrightarrow M_{k} \longrightarrow E_{k-1} \xrightarrow{\;\;
\partial_{k-1} \;\;} \dotsb \xrightarrow{\;\; \partial_{1} \;\;} E_{0}
\xrightarrow{\;\; \partial_{0} \;\;} M \xrightarrow{\;\; \partial_{-1}
\;\;} 0
\end{equation}
such that $E_{i} = \bigoplus_{1 \leq j \leq h_{i}} S(- \bm{m} - i
\bm{c})$, $M_{k} = \Ker \partial_{k-1}$ and $\bm{m} + k \bm{c} \in
\reg(M_{k})$.  Since $M = M_{0}$ and $M$ is $\bm{m}$-regular, the
first step in the induction holds.

Suppose $k \geq 0$.  Since $\mathbb{N} \scr{C} = \scr{K} =
\scr{K}^{\sat}$, Lemma~\ref{l:trunreg} and Lemma~\ref{l:trunsheaf}
show that $\reg(M_{k}) \subseteq \reg(M_{k}|_{(\bm{m} + k \bm{c} +
\mathbb{N} \scr{C})})$ and that $M_{k} / M_{k}|_{(\bm{m} + k\bm{c} +
\mathbb{N} \scr{C})}$ is a $B$-torsion module.  The induction
hypothesis states that $M_{k}$ is $(\bm{m} + k \bm{c})$-regular and
Theorem~\ref{t:gens} implies that $M_{k}|_{(\bm{m} + k \bm{c} +
\mathbb{N} \scr{C})}$ is generated in degree $\bm{m} + k \bm{c}$.  It
follows that there exists a surjective map $\partial_{k} \colon E_{k}
\longrightarrow M_{k}|_{(\bm{m} + k \bm{c} + \mathbb{N} \scr{C})}$
where $E_{k} := \bigoplus_{1 \leq j \leq b_{k}} S(- \bm{m} - k
\bm{c})$.  Setting $M_{k+1} := \Ker \partial_{k}$, we obtain the short
exact sequence
\begin{equation} \label{existscomplex:2}
0 \longrightarrow M_{k+1} \longrightarrow E_{k} \longrightarrow
M_{k}|_{(\bm{m} + k \bm{c} + \mathbb{N} \scr{C})} \longrightarrow 0 \,
.
\end{equation}
Combining \eqref{existscomplex:1} and \eqref{existscomplex:2} gives
the chain complex with $B$-torsion homology
\[
0 \longrightarrow M_{k+1} \longrightarrow E_{k} \xrightarrow{\;\;
\partial_{k} \;\;} \dotsb \xrightarrow{\;\; \partial_{1} \;\;} E_{0}
\xrightarrow{\;\; \partial_{0} \;\;} M \xrightarrow{\;\; \partial_{-1}
\;\;} 0
\]
where $E_{i} = \bigoplus_{1 \leq j \leq h_{i}} S(- \bm{m} - i \bm{c})$
and $M_{k+1} = \Ker \partial_{k}$.  It remains to show that $M_{k+1}$
is $\bigl( \bm{m} + (k+1) \bm{c} \bigr)$-regular.  Applying
Lemma~\ref{l:ses} to \eqref{existscomplex:2} yields
\[
\bigl( \bm{m} + k \bm{c} + \reg(S) \bigr) \cap \bigcap\limits_{1 \leq
i \leq \ell} \bigl( \bm{c}_{i} + \reg(M_{k}) \bigr) \subseteq
\reg(M_{k+1}) \, .
\]  
Our choice of $\bm{c}$ guarantees that $\bm{m} + (k+1)\bm{c}$ lies in
both $\bigcap\nolimits_{1 \leq i \leq \ell} \bigl( \bm{c}_{i} +
\reg(M_{k}) \bigr)$ and $\bm{m} + k \bm{c} + \reg(S)$ which completes
the proof.
\end{proof}

\begin{proof}[Proof of Theorem~\ref{i:complexes}]
In the introduction, $X$ was smooth so $\scr{K} = \scr{K}^{\sat}$.  We
also assumed that $\scr{C}$ was the set of minimal generators for
$\scr{K}$.  Hence, Theorem~\ref{t:seqtoreg} establishes Part~1 and
Theorem~\ref{t:existscomplex} proves Part~2.
\end{proof}

This theorem leads to a ``linear'' resolution of
$\mathscr{O}_{X}$-modules.

\begin{corollary}
Assume that $\mathbb{N} \scr{C} = \scr{K} = \scr{K}^{\sat}$.  If
$\bm{c} \in \reg(\mathscr{O}_{X}) \cap \bigl( \bigcap_{1 \leq j \leq
\ell} (\bm{c}_{j} + \scr{K}) \bigr)$ and $\bm{m} \in
\reg(\mathscr{F})$ then there is an exact sequence
\[
\;\; \dotsb \longrightarrow \mathscr{E}_{3} \longrightarrow
\mathscr{E}_{2} \longrightarrow \mathscr{E}_{1} \longrightarrow
\mathscr{E}_{0} \longrightarrow \mathscr{F} \longrightarrow 0
\] 
where $\mathscr{E}_{i} = \bigoplus\nolimits_{1 \leq j \leq h_{i}}
\mathscr{O}_{X}(- \bm{m} - i \bm{c})$.
\end{corollary}

\begin{proof}
Consider the $G$-graded $S$-module
\[
M = \Bigl( \bigoplus\nolimits_{\bm{p} \in G} H^{0} \bigl( X,
\mathscr{F}(\bm{p}) \bigr) \Bigr) \Bigr|_{(\bm{m} + \scr{K})} \, .
\]
Lemma~\ref{l:trunsheaf} implies that $\widetilde{M} = \mathscr{F}$ and
Proposition~\ref{p:sheafandmodulereg} shows that $M$ is
$\bm{m}$-regular.  Theorem~\ref{t:existscomplex} produces a chain
complex $\dotsb \longrightarrow E_{3} \longrightarrow E_{2}
\longrightarrow E_{1} \longrightarrow E_{0} \longrightarrow M
\longrightarrow 0$ with $B$-torsion homology and $E_{i} = \bigoplus_{1
\leq j \leq h_{i}} S(- \bm{m} - i \bm{c})$.  Applying the functor $F
\mapsto \widetilde{F}$, we obtain a chain complex of the desired form.
Moreover, Proposition~\ref{p:zerosheaf} implies that the homology of
this complex consists of the zero sheaf.  Therefore, the chain complex
is a locally free resolution of $\mathscr{F}$.
\end{proof}

We illustrate Theorem~\ref{t:existscomplex} with the following
examples.  The first example shows that Theorem~\ref{t:existscomplex}
is a converse to Theorem~\ref{t:seqtoreg} in the standard graded case.

\begin{example}
Let $S$ have the standard grading, let $B=\langle x_1, \dots,x_n
\rangle$, and let $\scr{C} = \{ \bm{1} \}$.  Since $\reg(S) =
\mathbb{N}$, we have $\bm{c} = \bm{1} \in \reg(S) \cap \bigl( \bm{1} +
\scr{K} \bigr)$.  Hence, if $\bm{m} \in \reg(M)$, then
Theorem~\ref{t:existscomplex} implies that there is a chain complex
$\mathbf{E} := \dotsb \longrightarrow E_{2} \longrightarrow E_{1}
\longrightarrow E_{0}$ with $B$-torsion homology and $E_{i} =
\bigoplus_{1 \leq j \leq h_{i}} S(- m - i)$ and there is a surjective
map $E_{0} \longrightarrow M|_{(\bm{m} + \mathbb{N})}$.  Conversely,
given a such chain complex, Example~\ref{e:standardsyzygies} shows
that $\bm{m} \in \reg(M)$.  In fact, Exercise~1.4.24 in \cite{BH}
proves that $\mathbf{E}$ must be a resolution.
\end{example}

\begin{example}
Let $M$ be the module in Example~\ref{e:basicRes} and let 
\[
\bm{c} := \bm{c}_{1} + \bm{c}_{2} + \bm{c}_{3} \in \reg(S) \cap \bigl(
\bigcap\nolimits_{1 \leq j \leq 3} (\bm{c}_{j} + \mathbb{N} \scr{C})
\bigr) \, .
\]
If $\bm{m} := \bm{0} \in \reg(M)$, then Theorem~\ref{t:existscomplex}
produces a chain complex with $B$\nobreakdash-torsion homology of the
form
\[
\dotsb \longrightarrow \bigoplus_{i} S(- 3\bm{c}) \longrightarrow
\bigoplus_{i} S(- 2\bm{c}) \longrightarrow \bigoplus_{1 \leq i \leq 7}
S(-\bm{c}) \longrightarrow S \, .
\]
However, applying Theorem~\ref{t:seqtoreg} to this chain complex
yields
\[
\bigcup_{\begin{subarray}{c} p_{1}, p_{2}, p_{3} \in
\mathbb{N} \\ p_{1} + p_{2} + p_{3} = 8
\end{subarray}} 
\left( \left[
\begin{smallmatrix} p_{1} \\ p_{2} \\ p_{3} \end{smallmatrix} \right] +
\mathbb{N}^{3} \right) \subset \reg(M), \, 
\]
and the smaller set does not contain $\bm{0}$.
\end{example}


\def\cprime{$'$}
\providecommand{\bysame}{\leavevmode\hbox to3em{\hrulefill}\thinspace}

\end{document}